\documentclass[reqno,11pt,letterpaper]{amsart}
\usepackage{amsmath,amsfonts,amsthm,amscd,amssymb}
\usepackage[T1]{fontenc}
\usepackage{lmodern}
\usepackage{thmtools,thm-restate}
\usepackage[backend=biber,style=numeric]{biblatex}
\usepackage{indentfirst,bm,mathtools}
\usepackage{tikz}
\usepackage{booktabs,array,enumitem}
\usepackage{microtype,fancyhdr}
\usepackage[letterpaper, inner=30mm, outer=27mm, top=28mm, bottom=25mm]{geometry}
\usepackage[colorlinks=true,linkcolor=blue!45!black,citecolor=blue!45!black,urlcolor=blue!45!black]{hyperref}
\newtheorem{theorem}{Theorem}[section]
\newtheorem{proposition}{Proposition}[section]
\newtheorem{lemma}{Lemma}[section]
\newtheorem{corollary}{Corollary}[section]

\theoremstyle{definition}

\theoremstyle{remark}

\numberwithin{equation}{section}

\usetikzlibrary{arrows.meta,calc,positioning,decorations.markings,bending}
\definecolor{surfaceblue}{RGB}{91,153,194}
\definecolor{surfaceedge}{RGB}{30,79,113}
\definecolor{glueboundary}{RGB}{119,74,151}
\newcommand{\R}{\mathbb R}

\newcommand{\Z}{\mathbb Z}
\newcommand{\dd}{\,\mathrm d}
\newcommand{\Reg}{\operatorname{Reg}}
\newcommand{\Sing}{\operatorname{Sing}}
\newcommand{\supp}{\operatorname{supp}}
\newcommand{\diam}{\operatorname{diam}}
\newcommand{\dist}{\operatorname{dist}}
\newcommand{\vol}{\operatorname{vol}}
\newcommand{\rank}{\operatorname{rank}}

\newcommand{\cP}{\mathcal P}
\newcommand{\cQ}{\mathcal Q}
\newcommand{\cM}{\mathcal M_1(\phi)}
\newcommand{\eps}{\varepsilon}
\newcommand{\indic}{\mathbf1}
\newcommand{\mass}[1]{\lVert#1\rVert}
\begin{document}
\title[Measure of Maximal Entropy]{Measures of Maximal Entropy for Geodesic Flows on Noncompact Manifolds of Nonpositive Curvature}
\author{Fei Liu}
\address{College of Mathematics and Systems Science, Shandong University of Science and Technology, Qingdao, 266590, China}
\email{liufei@math.pku.edu.cn}
\author{Xiaokai Liu}
\address{Department of Mathematics, Southern University of Science and Technology, Shenzhen, 518055, China}
\email{liuxk2024@mail.sustech.edu.cn}
\thanks{Corresponding author: Xiaokai Liu.}
\date{\today}
\begin{abstract}
We establish the variational characterization of the measure of maximal entropy on noncompact manifolds of nonpositive sectional curvature. Assuming that the deck transformation group contains an axial isometry with a rank one axis and has a finite critical exponent, we prove that the supremum of the entropies of invariant probabilities equals the critical exponent and is attained if and only if the Bowen--Margulis--Sullivan (BMS) measure is finite. In that case, the measure of maximal entropy is unique and equals the normalized BMS measure.

The proof combines conditional-entropy rigidity on the regular set with an HTS-based entropy exclusion on the singular set, overcoming both non-compactness and the loss of uniform hyperbolicity due to zero curvature. We also obtain mixing under a non-arithmetic length-spectrum assumption.  For nonpositively curved surfaces of genus at least two, the existence of a closed rank one geodesic and the non-arithmeticity of the axial length spectrum are automatic.
\end{abstract}
\maketitle
\section{Introduction}\label{sec:introduction}
The relation between orbit growth and entropy is a central theme in the study of geodesic flows. On closed manifolds, the variational principle expresses topological entropy as the supremum of the measure-theoretic entropies of invariant probability measures. A probability attaining this supremum is called a measure of maximal entropy. Its existence and uniqueness link the flow's orbit complexity to a specific invariant probability measure.

Margulis constructed invariant measures from measures on stable and unstable leaves in~\autocite{Margulis1969}. Bowen approached the problem through entropy and periodic orbits, constructing invariant probabilities from periodic-orbit measures and establishing the uniqueness of the measure of maximal entropy for closed negatively curved manifolds~\cite{Bowen1972Periodic,Bowen1972Equidistribution}. A different construction arose from Patterson's study of limit sets of Fuchsian groups and Sullivan's extension to discrete groups of isometries of real hyperbolic space~\cite{Patterson1976,Sullivan1979}. Using these conformal boundary measures, Sullivan constructed an invariant measure for the geodesic flow on the unit tangent bundle of the quotient. For closed negatively curved manifolds, its normalization coincides with the Bowen--Margulis measure; see \cite[\S3.5]{Kaimanovich1990}. In this article, we refer to this measure as the Bowen--Margulis--Sullivan (BMS) measure. On noncompact manifolds, this measure may be infinite, raising the question of whether the finiteness of the BMS measure implies the existence of a measure of maximal entropy.

In nonpositive curvature, flat strips and parallel Jacobi fields prevent the geodesic flow from being uniformly hyperbolic. Manning identified topological entropy with volume growth for closed nonpositively curved manifolds \cite{Manning1979}, and Knieper proved the uniqueness of the measure of maximal entropy for closed rank one manifolds using the conformal density on the ideal boundary \cite{Knieper1998}. Ricks extended this uniqueness theorem to compact, geodesically complete, locally CAT(0) spaces admitting a weak rank one axis in the universal cover \cite{Ricks2021}. For closed rank one manifolds, Burns, Climenhaga, Fisher, and Thompson established uniqueness of equilibrium states for H\"older potentials whose pressure on the singular set is strictly smaller than the full pressure \cite{BCFT2018}. Recently, Wu constructed these equilibrium states from boundary measures, subject to additional hypotheses on the potential \cite{Wu2025}.

On a noncompact manifold, the invariant measures need not be finite, and the entropy supremum may not be attained. For complete manifolds with pinched negative curvature and non-elementary fundamental group, Otal and Peign\'e established the variational characterization: the supremum of the entropies of invariant probabilities equals the critical exponent and is attained if and only if the BMS measure is finite, in which case its normalization is the unique measure of maximal entropy \cite[Theorems~1 and~2]{OtalPeigne2004}. The assumption of uniformly bounded curvature derivatives used in their proof was omitted from the original statement; see the later clarification in \cite[footnotes~4 and~7]{BAPP2021}. Roblin developed the recurrence, mixing, and equidistribution theory for discrete groups acting on proper CAT($-1$) spaces \cite{Roblin2003}. In particular, his work relates divergence of the Poincar\'e series at $\delta_{\Gamma}$ to conservativity and ergodicity, and establishes mixing and orbital asymptotics under a non-arithmetic length-spectrum hypothesis. These results supply the geometric and ergodic background for the noncompact variational problem.

The finiteness of the BMS measure has also been studied. Yue obtained finiteness criteria in pinched negative curvature \cite{Yue1996}. Dal'bo, Otal, and Peign\'e obtained a parabolic summability criterion for geometrically finite divergent groups in pinched negative curvature \cite{DOP2000}. For H\"older potentials of finite pressure in pinched negative curvature, Pit and Schapira studied finiteness of Gibbs measures through recurrence and summability of excursions outside compact sets \cite{PitSchapira2018}. For pinched negatively curved manifolds with bounded first derivative of curvature, Gou\"ezel, Schapira, and Tapie related pressure at infinity to strong positive recurrence and proved that a strict pressure gap at infinity for a H\"older potential of finite pressure implies finiteness of its Gibbs measure \cite{GST2023}. Beyond negative curvature, Liu, Liu, and Wang obtained a finiteness criterion on noncompact manifolds without conjugate points satisfying uniform visibility and Axiom~2 \cite{LLW2025Finiteness}. 

For noncompact manifolds of nonpositive curvature, Link and Picaud established the Hopf--Tsuji--Sullivan (HTS) dichotomy for non-elementary groups containing a rank one isometry \cite{LP16}. In the divergent case, the BMS measure is conservative and ergodic and gives full measure to the rank one set. Link extended this dichotomy to Hadamard spaces with a weak rank one group action and compared the Ricks and Knieper measures in the conservative manifold case \cite{Link2018}. For complete uniform visibility manifolds without conjugate points satisfying Axiom~2 and the non-elementarity assumption, Liu, Liu, and Wang established an HTS dichotomy \cite{LLW2025HTS}. Liu and Wang studied uniqueness, up to scale, within the class of Knieper measures when a finite such measure exists \cite{LiuWang2021}. Uniqueness within this geometrically defined class does not establish maximal entropy or uniqueness among all invariant probabilities. The latter comparison must also account for probabilities supported on the singular set.

Mixing in this setting already has a general theory. Babillot developed a stable and unstable leaf argument for mixing \cite{Babillot02}, and Ricks related mixing to the length spectrum in CAT(0) spaces admitting a weak rank one isometry \cite{Ricks17}. Later, Link proved mixing under a non-arithmetic length-spectrum hypothesis for proper, geodesically complete weak rank one Hadamard spaces \cite{Link2020}. In the manifold setting with finite BMS measure, Burniol Clotet explicitly related mixing to non-arithmeticity of the rank one length spectrum \cite[Theorem 2.1]{BurniolClotet2023}. In the compact setting beyond nonpositive curvature, Liu, Liu, and Wang proved mixing of the unique measure of maximal entropy for smooth rank one manifolds without focal points, and the Bernoulli property for compact surfaces of genus greater than one without focal points \cite{LLW2021Mixing}. The mixing statement below follows from this theory after identifying the BMS measure with Ricks' measure.

The main result of this paper extends the Otal--Peign\'e variational characterization to complete manifolds of nonpositive curvature under a closed rank one geodesic hypothesis and finiteness of the critical exponent. The finiteness of the critical exponent is automatic in pinched negative curvature. The assumptions allow non-compactness, infinite volume, infinitely generated fundamental groups, and sectional curvature unbounded below. We identify the supremum of the entropies of invariant probabilities with the critical exponent and show that every maximizing probability is the normalized BMS measure. 

The extension requires simultaneous control of non-compactness and the degeneracy caused by zero curvature. For regular invariant probabilities, we adapt the conditional-entropy method of Otal and Peign\'e \cite[\S\S2--3]{OtalPeigne2004} by replacing uniform contraction with contraction accumulated during returns to compact windows where the unstable Green operator is positive definite. Finite-entropy return and deck coding, together with a flow-coordinate refinement, yield the required unstable partitions. This establishes entropy rigidity without assuming finiteness of the BMS measure: a regular maximizing probability forces finiteness and coincides with the normalized BMS measure.

The singular set presents an additional obstruction absent in negative curvature. We use the rank one HTS dichotomy \cite{LP16} to convert conservative recurrence into a strict entropy deficit for every invariant probability concentrated on this set in the divergent case. The entropy-deficit comparison follows the compact arguments of Knieper and Ricks \cite{Knieper1998,Ricks2021}, but here the candidate probability may have noncompact support and the BMS measure may be infinite. We fix a positive-mass compact return set and increase the sampling time. Geodesics tracking long singular segments avoid a fixed regular window, whereas typical BMS orbits return to that window. Combined with the overlap estimates of Link and Picaud \cite[Section~6]{LP16}, this gives a vanishing factor in the return-label count that is uniform in the return step. The argument requires no uniform entropy gap on the singular set and applies even to infinite conservative BMS measures, completing the characterization of maximal entropy measures by finiteness of the BMS measure.

We state our main results here:

\begin{theorem}\label{thm:main}
    Let $M=X/\Gamma$ be a connected, smooth, complete Riemannian manifold of dimension $n\ge2$ with sectional curvature $K\le0$, where $X$ is its universal cover and $\Gamma$ is its deck transformation group. Assume that $\Gamma$ is non-elementary with the critical exponent $0<\delta_{\Gamma}<\infty$, and $M$ contains a closed rank one geodesic. Here rank one means that a lift of the geodesic to $X$ admits no nonzero parallel normal Jacobi field. Fix a Patterson--Sullivan density of dimension $\delta_{\Gamma}$, and let $m$ be the associated Bowen--Margulis--Sullivan measure.

Let $\cM$ denote the set of probability measures invariant under the entire geodesic flow $\phi$, and let $h_{\mathrm{top}}^{\,c}(\phi)$ denote the supremum of the topological entropies of its restrictions to compact flow-invariant subsets. Then
\begin{displaymath}
    h_{\mathrm{top}}^{\,c}(\phi)=H:=\sup_{\nu\in\cM}h_\nu(\phi^1)=\delta_\Gamma.
\end{displaymath}
A measure of maximal entropy exists if and only if $m$ is finite. When it exists, it is unique and equals $\bar m=m/\mass m$. In particular,
\begin{displaymath}
    h_{\bar m}(\phi^1)=\delta_\Gamma.
\end{displaymath}
\end{theorem}

The theorem requires neither finite Riemannian volume nor a global lower bound on sectional curvature, and the fundamental group need not be finitely generated. Thus the variational characterization extends to noncompact manifolds with flat directions and curvature unbounded below. The proof also excludes singular maximizing probabilities without assuming a uniform entropy gap on the singular set; in the divergent case, this exclusion remains valid even when the BMS measure is infinite.

For an axial element $\gamma\in\Gamma$, let $\ell(\gamma)=\inf_x d(x,\gamma x)>0$ be its translation length. The axial length spectrum is
\begin{displaymath}
    \mathcal L(\Gamma)=\{\ell(\gamma):\gamma\in\Gamma~\text{is axial}\}.
\end{displaymath}
It is non-arithmetic if it is not contained in $a\Z$ for any $a>0$, or equivalently if the closed additive subgroup it generates is $\R$.

\begin{theorem}\label{thm:scope:mixing}
    Under the assumptions of Theorem~\ref{thm:main}, suppose that $m$ is finite and $\mathcal L(\Gamma)$ is non-arithmetic. Then the geodesic flow is mixing with respect to $\bar m$:
    \begin{displaymath}
        \bar m(A\cap \phi^{-t}B)\longrightarrow\bar m(A)\bar m(B)\qquad(t\to+\infty)
    \end{displaymath}
    for all Borel sets $A,B\subset T^1M$.
\end{theorem}

For orientable surfaces of genus at least two, the closed rank one hypothesis and non-arithmeticity of the axial length spectrum are automatic. This follows from the surface results of Link, Peign\'e, and Picaud~\cite{LPP2006,LP16}, applied to a finite-type covering surface; see Section~\ref{sec:surfaces}. We therefore obtain the following corollary.

\begin{corollary}\label{cor:scope:surface}
    Let $S$ be a connected, orientable, smooth, complete surface without boundary, of genus at least two and with $K\le0$. Suppose that $\delta_{\pi_1(S)}<\infty$. Then
    \begin{displaymath}
        h_{\mathrm{top}}^{\,c}(\phi)=\sup_{\nu\in\cM}h_\nu(\phi^1)=\delta_{\pi_1(S)}>0.
    \end{displaymath}
    The geodesic flow has a measure of maximal entropy if and only if its BMS measure is finite. In that case, the measure is unique, equals the normalized BMS measure, and is mixing.
\end{corollary}

The article is organized as follows: Section~\ref{sec:preliminaries} fixes the geometric and entropy conventions and states the recurrence results used below. Section~\ref{sec:entropy} identifies the compact topological entropy with the critical exponent. Section~\ref{sec:uniqueness} proves the finite BMS entropy formula, uniqueness, and the finiteness criterion. Finally, Section~\ref{sec:mixing} proves mixing for finite BMS measures and derives the surface corollary.

\section{Setting and preliminaries}\label{sec:preliminaries}\label{sec:geometry}

Let $M=X/\Gamma$ be a connected, complete, smooth Riemannian manifold of dimension $n\ge2$ with sectional curvature $K\le0$. The universal cover $X$ is a Hadamard manifold, and the deck transformation group $\Gamma$ acts freely and properly discontinuously by isometries. Unless a weaker hypothesis is explicitly stated, we assume that $\Gamma$ is non-elementary, contains an axial isometry with a rank one axis, and has a finite positive critical exponent. Non-elementary means that the limit set is infinite and that $\Gamma$ fixes no point of the ideal boundary. We write $\pi:T^1 X\to X$ for the standard projection and use the same notation on the quotient. Denote the geodesic flow by $\phi^tv=c_v'(t)$, where $c_v$ is the unit-speed geodesic with initial vector $v$.

\subsection{Geometry and rank one}\label{sec:geometry:rank}
The Cartan--Hadamard theorem identifies $X$ with each tangent space by the exponential map. Closed bounded subsets of $X$ are compact. The distance between two geodesics parametrized on the same interval is convex. We use these standard facts in the form given in \cite[Chapter~I, Sections~2~and~4,~Proposition~5.4]{Ballmann95}; they do not require a lower curvature bound or a compact quotient.

The ideal boundary $\partial X$ consists of asymptotic classes of geodesic rays. For each $o\in X$, a boundary point is represented by a unique ray issuing from $o$. The cone topology identifies $\partial X$ with $S^{n-1}$ and $X\cup\partial X$ with a closed ball. We denote the endpoints of a complete geodesic $c$ by $c^-$ and $c^+$. The limit set $\Lambda_\Gamma$ is the set of accumulation points of $\Gamma o$ in $\partial X$ under the cone topology.

A normal Jacobi field along $c$ satisfies
\begin{displaymath}
    J''+R(J,c')c'=0,\qquad J\perp c'.
\end{displaymath}
The rank of the unit tangent vector along $c$ is defined as one plus the dimension of its space of parallel normal Jacobi fields. Set
\begin{displaymath}
    \Sing=\{v\in T^1M:\rank(v)>1\},\qquad\Reg=T^1M\setminus\Sing.
\end{displaymath}
The rank function is upper semicontinuous and invariant under the geodesic flow; see \cite[Chapter~IV, Section~4]{Ballmann95}. Thus $\Reg$ is open, and $\Sing$ is closed, and both sets are flow invariant. A geodesic is \emph{weakly rank one} if it does not bound a flat half-plane. Rank one implies weak rank one. Moreover, the flat strip theorem \cite[Chapter~I, Corollary~5.8(ii)]{Ballmann95} shows that a rank one endpoint pair determines a unique unparametrized geodesic: a strip between distinct parallel lines would supply a nonzero parallel normal Jacobi field.

A rank one geodesic has endpoint neighborhoods such that every pair of points in these neighborhoods is joined by a rank one geodesic meeting a fixed compact ball~\cite[Lemma 2.1]{LP16}. After fixing a time origin by a Busemann level, the connecting geodesics depend continuously on their endpoints: the compact intersection bound gives convergence, and uniqueness identifies every limit. A weak rank one axial isometry has north--south dynamics on the ideal boundary \cite[Lemma~2.4]{LP16}. In the non-elementary case, there are two such elements with disjoint fixed-point sets \cite[Section~3]{LP16}; sufficiently large powers generate a free semigroup by ping-pong.

\subsection{The critical exponent and entropy}
For $o\in X$, put
\begin{displaymath}
    \delta_\Gamma=\limsup_{R\to\infty}\frac1R\log\#\{\gamma\in\Gamma:d(o,\gamma o)\le R\}.
\end{displaymath}
We write $N_o(R)=\#\{\gamma\in\Gamma:d(o,\gamma o)\le R\}$ for the orbital counting function. The inequality
\begin{displaymath}
    |d(o,\gamma o)-d(o',\gamma o')|\le2d(o,o')
\end{displaymath}
shows that $\delta_{\Gamma}$ is independent of the choice of $o$. Proper discontinuity ensures that every finite-radius count is finite.

For every $s>\delta_{\Gamma}$ and every fixed $D\geq0$, the definition of the critical exponent gives a constant $C_{s,D}\geq1$ such that
\begin{equation}\label{eq:fb:orbitalcount}
    N_o(r+2D)\leq C_{s,D}e^{sr}\qquad(r\geq0).
\end{equation}
Increasing the constant covers the bounded interval outside the asymptotic estimate.

The Poincar\'e series
\begin{displaymath}
    \mathcal P_\Gamma(s)=\sum_{\gamma\in\Gamma}e^{-s d(o,\gamma o)}
\end{displaymath}
converges for $s>\delta_{\Gamma}$ and diverges for $s<\delta_{\Gamma}$. We say that $\Gamma$ is of convergent or divergent type according to the behavior of this series at $s=\delta_{\Gamma}$. If $a,b\in\Gamma$ freely generate a positive semigroup and $L=\max\{d(o,ao),d(o,bo)\}$, the $2^k$ positive words of length $k$ have displacement at most $kL$, giving
\begin{equation}\label{eq:prelim:free_growth}
    \delta_{\Gamma}\ge\frac{\log2}{L}>0.
\end{equation}
Thus $\delta_{\Gamma}>0$ whenever $\Gamma$ is non-elementary and contains a weak rank one axial isometry.

Let $\cM$ denote the set of flow-invariant Borel probability measures, and define the probability variational entropy by
\begin{displaymath}
    H=\sup_{\nu\in\cM}h_\nu(\phi^1).
\end{displaymath}
A measure of maximal entropy is a member of $\cM$ attaining this supremum. To specify the topological quantity on the noncompact phase space, define
\begin{equation}\label{eq:prelim:compact_entropy}
    h_{\mathrm{top}}^{\,c}(\phi)=\sup_{\substack{C\subset T^1M~\text{compact}\\\phi^tC=C\ (t\in\R)}}h_{\mathrm{top}}(\phi^1|_C).
\end{equation}
On a compact invariant set $C$, separated orbit segments define the entropy on the right: if $s_C(T,\varepsilon)$ is the largest cardinality of a set separated for $d_T(v,w)=\max_{0\le t\le T}d(\phi^tv,\phi^tw)$, then
\begin{displaymath}
    h_{\mathrm{top}}(\phi^1|_C)=\lim_{\varepsilon\downarrow0}\limsup_{T\to\infty}\frac1T\log s_C(T,\varepsilon).
\end{displaymath}
This value is independent of the choice of a compatible metric on $C$. The compact variational principle identifies it with the supremum of $h_\nu(\phi^1)$ over $\phi^1$-invariant probabilities supported on $C$ \cite[Theorem~3.4.1]{PU2010}. Averaging over one unit of flow time preserves entropy and gives a flow-invariant probability by~\eqref{eq:lower:time_average}, so $h_{\mathrm{top}}^{\,c}(\phi)\le H$. For compact $M$, definition~\eqref{eq:prelim:compact_entropy} gives the usual topological entropy.

\subsection{BMS measure and the Hopf--Tsuji--Sullivan dichotomy}
For a ray $c$ ending at $\xi$, use the Busemann convention
\begin{displaymath}
    \beta_\xi(x,y)=\lim_{t\to\infty}\left(d(x,c(t))-d(y,c(t))\right).
\end{displaymath}
It satisfies the standard cocycle, distance, and equivariance identities \cite[Section~2, equations~(3)--(4)]{LP16}
\begin{displaymath}
    \begin{aligned}
        \beta_\xi(x,z)&=\beta_\xi(x,y)+\beta_\xi(y,z),\\
        |\beta_\xi(x,y)|&\le d(x,y),\\
        \beta_{\gamma\xi}(\gamma x,\gamma y)&=\beta_\xi(x,y).
    \end{aligned}
\end{displaymath}
The Patterson construction gives a $\delta_{\Gamma}$-dimensional conformal density $(\sigma_x)_{x\in X}$ with
\begin{displaymath}
    \gamma_*\sigma_x=\sigma_{\gamma x},\qquad\frac{\dd\sigma_x}{\dd\sigma_y}(\xi)=e^{-\delta_{\Gamma}\beta_\xi(x,y)},\qquad \supp\sigma_x=\Lambda_\Gamma.
\end{displaymath}
We fix this density and normalize it by $\sigma_o(\partial X)=1$; existence and full support in the present setting are given in \cite[Section 4 and Lemma 4.1]{LP16}.

Write $\partial^2X$ for the ordered endpoint pairs joined by complete geodesics. For $(\xi,\eta)\in\partial^2X$, let $P(\xi,\eta)$ be the union of these geodesics and set
\begin{displaymath}
    G_o(\xi,\eta)=\beta_\xi(o,z)+\beta_\eta(o,z),\qquad z\in P(\xi,\eta).
\end{displaymath}
This is twice the Gromov product in \cite[Section~2, equation~(5)]{LP16}, so it is independent of $z$ and nonnegative. If $P(\xi,\eta)$ meets $B(o,R)$, the Busemann distance bound gives $0\le G_o\le2R$.

The parallel-set theorem \cite[Chapter~I, Proposition~5.9]{Ballmann95} gives the splitting $P(\xi,\eta)=C(\xi,\eta)\times\R$. It carries the Riemannian volume of its intrinsic dimension, with a convex boundary allowed. This dimension can depend on the endpoint pair. Following the geometric construction in \cite[Section 2 and Section 6, equation (22)]{LP16}, define the measure on $T^1X$ as
\begin{equation}\label{eq:full_knieper}
    \widetilde m(E)=\int_{\partial^2X} e^{\delta_{\Gamma} G_o(\xi,\eta)}\vol_{P(\xi,\eta)}\{z:v_{z,\eta}\in E\}\dd\sigma_o(\xi)\dd\sigma_o(\eta),
\end{equation}
where $v_{z,\eta}$ is the vector based at $z$ pointing toward $\eta$. The fiber measurability and the invariance of $\widetilde m$ under $\Gamma$ and the geodesic flow follow from \cite[Section~6, equation~(22)]{LP16}. Full support of the conformal density and the endpoint neighborhoods of a rank one axis gives $\widetilde m\ne0$. Its quotient is the BMS measure, denoted by $m$.

The measure is locally finite. Fix $R>0$ and, for each parallel set $P$ meeting $B(o,R)$, choose a relative interior point $z\in P\cap B(o,R+1)$. Convexity gives
\begin{displaymath}
    P\cap B(o,R)\subset
    \exp_z\{u\in T_zP:|u|\le2R+1\}.
\end{displaymath}
The differential of the exponential map is uniformly bounded on the compact set $\{(z,u):z\in\overline{B(o,R+1)},\ |u|\le2R+1\}$. Since $1\le\dim P\le n$, the area formula gives a constant $C_R$ independent of $P$ with $\vol_P(P\cap B(o,R))\le C_R$. Integrating against the normalized endpoint measures yields
\begin{equation}\label{eq:prelim:local_finiteness}
    \widetilde m(T^1B(o,R))\le C_Re^{2\delta_{\Gamma} R}<\infty.
\end{equation}
Covering charts give local finiteness on $T^1M$, hence a Radon, $\sigma$-finite measure.

For a regular vector put $\xi=v^-$, $\eta=v^+$, and $t=\beta_\xi(\pi v,o)$. These Hopf coordinates satisfy $t(\phi^sv)=t(v)+s$, and uniqueness of the connecting line reduces~\eqref{eq:full_knieper} to
\begin{equation}\label{eq:hopf_bms}
    \dd\widetilde m=e^{\delta_{\Gamma} G_o(\xi,\eta)}
    \dd\sigma_o(\xi)\dd\sigma_o(\eta)\dd t.
\end{equation}
We write $\mass m=m(T^1M)$ and $\bar m=m/\mass m$ when it is finite.

A point $\xi\in\Lambda_\Gamma$ is radial (or conical) if a ray ending at $\xi$ returns to a fixed-radius neighborhood of $\Gamma o$ at arbitrarily large positive times. Let $\Lambda_{\mathrm{rad}}$ be the radial limit set.

\begin{theorem}[Link--Picaud]\label{thm:geo:HTS}
    Under the standing hypotheses, exactly one of the following alternatives holds:
    \begin{enumerate}
        \item $\mathcal P_\Gamma(\delta_{\Gamma})<\infty$, $\sigma_o(\Lambda_{\mathrm{rad}})=0$, and the geodesic flow is completely dissipative with respect to $m$;
        \item $\mathcal P_\Gamma(\delta_{\Gamma})=\infty$, $\sigma_o(\Lambda_{\mathrm{rad}})=1$, and the flow is conservative and ergodic with respect to $m$. Moreover, $m(\Sing)=0$.
    \end{enumerate}
\end{theorem}
This is \cite[Theorem 7.2 and Corollary 1]{LP16}. Conservativity and ergodicity here concern the $\sigma$-finite measure $m$. In the dissipative alternative, for every nonnegative $\varphi\in L^1(m)$, the integral $\int_0^\infty\varphi(\phi^tv)\dd t$ is finite for $m$-almost every $v$.

We record the consequences of this dichotomy in the forms used below. Recurrence for a nonzero finite invariant measure excludes complete dissipativity, so Theorem~\ref{thm:geo:HTS} gives
\begin{equation}\label{eq:prelim:finite_BMS}
    \mass m<\infty\quad\Longrightarrow\quad\Gamma~\text{divergent},~\bar m~\text{ergodic},~\bar m(\Reg)=1.
\end{equation}
In the divergent case, conservativity and ergodicity imply, for every open $U$ with $0<m(U)<\infty$,
\begin{equation}\label{eq:prelim:HTS_visits}
    m\{v:\{t>0:\phi^tv\in U\}\text{ is bounded}\}=0.
\end{equation}
Indeed, almost every point of $U$ returns infinitely often, and ergodicity makes the flow saturation of $U$ co-null; openness allows both assertions to be expressed in terms of rational times.

We shall also use the standard invariant-density consequence
\begin{equation}\label{eq:prelim:ac_BMS}
    \nu\in\cM,\quad \nu\ll m
    \quad\Longrightarrow\quad \mass m<\infty\text{ and }\nu=\bar m.
\end{equation}
To see why the dissipative alternative is excluded here, choose a bounded everywhere positive $\varphi\in L^1(m)$, using $\sigma$-finiteness. Complete dissipativity makes its forward integral finite $m$-almost everywhere, hence $\nu$-almost everywhere, so its time average tends to zero. The probability ergodic theorem for $\nu$ contradicts this because the limiting conditional expectation has a positive integral. Thus $m$ is conservative and ergodic. Invariance of $\nu=fm$ makes $f$ flow invariant modulo $m$, so ergodicity makes $f$ a positive constant and normalization gives~\eqref{eq:prelim:ac_BMS}.

\subsection{Entropy conventions and standard formulas}\label{sec:entropy_tools}
We assume that all probability spaces below are standard Borel, and all transformations are invertible and probability preserving. In the article, $M$ and $T^{1}M$ are standard Borel spaces. We use Rokhlin disintegration, the ergodic theorems, and the entropy identities in \cite{Walters82,PU2010,AbramovInd59,Abramov59}. Completions of the standard Borel probability spaces are Lebesgue spaces, as required in the partition formulation of \cite[Sections~2.6--2.8]{PU2010}.

For a countable Borel partition $\mathcal P=\{P_i\}$, set
\begin{displaymath}
    H_\nu(\mathcal P)=-\sum_i\nu(P_i)\log\nu(P_i).
\end{displaymath}
Here we adopt the convention $0\log0=0$. Given a standard Borel probability space $(Y,\nu)$ and a Borel map $c:Y\to Z$ into a standard Borel space $Z$, we also consider the fiber partition $\alpha=\{c^{-1}(z):z\in Z\}$. This partition represents the information $\sigma(\alpha)=c^{-1}\mathcal B(Z)$. For two partitions $\alpha$ and $\beta$, we write $\alpha\succeq\beta$ when $\alpha$ refines $\beta$. For $\lambda=c_*\nu$, disintegration gives probability kernels $\nu_z$ such that
\begin{displaymath}
    \nu(E)=\int_Z\nu_z(E)\dd\lambda(z),\qquad
    \nu_z(c^{-1}(z))=1\quad\text{for almost every}~z.
\end{displaymath}
Conditional probabilities are identified modulo $\nu$-null sets. The underlying Borel fibers are retained when geometric reference measures are placed on them.

For a sub-$\sigma$-algebra $\mathcal A$ and countable $\mathcal P$, conditional entropy is the integral of $-\sum_i p_i\log p_i$, where $p_i=\mathbb E_\nu(\indic_{P_i}\mid\mathcal A)$. For fiber partitions, define
\begin{displaymath}
    H_\nu(\alpha\mid\beta)=
    \sup_{\substack{\mathcal E~\text{finite}\\\sigma(\mathcal E)\subset\sigma(\alpha)}}H_\nu(\mathcal E\mid\sigma(\beta)).
\end{displaymath}
Any increasing sequence of finite partitions generating $\sigma(\alpha)$ realizes this supremum as a limit. On a conditional fiber, a countably supported distribution has its usual Shannon entropy.

Conditional entropy satisfies the chain rule
\begin{displaymath}
    H_\nu(\alpha\vee\gamma\mid\beta)=H_\nu(\alpha\mid\beta)+H_\nu(\gamma\mid\alpha\vee\beta), 
\end{displaymath}
and decreases as the conditioning information increases. If $\alpha_j$ increases to $\alpha$, then $H_\nu(\alpha_j\mid\beta)$ increases to $H_\nu(\alpha\mid\beta)$. For $\mathcal A_j\uparrow\mathcal A$ and $H_\nu(\mathcal P)<\infty$, the standard continuity theorem gives
\begin{equation}\label{et:conditioning_limit}
    H_\nu(\mathcal P\mid\mathcal A_j)\downarrow H_\nu(\mathcal P\mid\mathcal A).
\end{equation}
These formulas are understood in $[0,\infty]$ and follow by finite approximation from the standard entropy identities; see \cite[Theorem~2.3.1 and Theorems~2.8.5--2.8.6]{PU2010}.

Let $f^k\mathcal P=\{f^kP:P\in\mathcal P\}$. This partition records the symbol of $f^{-k}x$. Put
\begin{displaymath}
    \mathcal P^{(n)}=\bigvee_{j=0}^{n-1}f^{-j}\mathcal P,\qquad \mathcal P^- =\bigvee_{k=0}^\infty f^k\mathcal P.
\end{displaymath}
For $H_\nu(\mathcal P)<\infty$, the chain rule and~\eqref{et:conditioning_limit} give
\begin{align*}
    h_\nu(f,\mathcal P)=\lim_{n\to\infty}\frac1nH_\nu(\mathcal P^{(n)})=H_\nu(f^{-1}\mathcal P\mid\mathcal P^-),\\
    H_\nu\left(\bigvee_{j=1}^nf^{-j}\mathcal P\,\middle|\,\mathcal P^-\right)=n h_\nu(f,\mathcal P).
\end{align*}
The Kolmogorov--Sinai entropy $h_\nu(f)$ is the supremum over finite Borel partitions. Finite coarse approximations show that countable finite-entropy partitions also have entropy rate at most $h_\nu(f)$.

Let $N$ be a positive integer-valued random variable with finite mean, and suppose that a label $J$ has at most $Ce^{an}$ possible values conditional on $N=n$, where $C\ge1$ and $a\ge0$. Gibbs' inequality with the reference distribution $q_n=2^{-n}$ gives $H(N)\le(\log 2)\mathbb E N$, while $H(J\mid N=n)\le\log C+an$. Thus, the chain rule gives
\begin{equation}\label{et:finite_labels}
    H(N,J)=H(N)+H(J\mid N)\le(\log 2+a)\mathbb E N+\log C.
\end{equation}

For an invariant probability $\nu$, let $\mathcal I_f$ and $\mathcal I_\phi$ denote the invariant $\sigma$-algebras of a time map and of the full flow, respectively. Recurrence gives infinitely many forward and backward returns to a positive-measure set for almost every point of that set. For $\varphi\in L^1(\nu)$, the ergodic theorems give
\begin{displaymath}
    \frac1n\sum_{j=0}^{n-1}\varphi(f^jx)\longrightarrow\mathbb E_\nu(\varphi\mid\mathcal I_f)(x),\qquad\frac1T\int_0^T\varphi(\phi^tx)\dd t\longrightarrow\mathbb E_\nu(\varphi\mid\mathcal I_\phi)(x)
\end{displaymath}
almost everywhere. These limits are constant when the corresponding system is ergodic; see \cite{Walters82}.

Let $A$ have measure $p>0$, let $B=\bigcup_{j\in\Z}f^jA$ have measure $s$, and write $\nu_A=\nu|_A/p$ and $\nu_B=\nu|_B/s$. On the recurrent part of $A$, define
\begin{displaymath}
    N(x)=\min\{n\ge1:f^nx\in A\},\qquad F(x)=f^{N(x)}x.
\end{displaymath}
The return map preserves $\nu_A$. The return tower partitions $B$ modulo null sets and gives, for every nonnegative measurable $\psi$,
\begin{displaymath}
    \int_A\sum_{j=0}^{N(x)-1}\psi(f^jx)\dd\nu(x)=\int_B\psi\dd\nu.
\end{displaymath}
In particular, Kac's formula is
\begin{equation}\label{et:kac}
    \int_A N\dd\nu=s,\qquad \int_A N\dd\nu_A=s/p.
\end{equation}
Applying Abramov's induced entropy formula to $(B,\nu_B,f)$ yields
\begin{equation}\label{et:abramov_induced}
    p\,h_{\nu_A}(F)=s\,h_{\nu_B}(f).
\end{equation}
The non-ergodic version follows by applying the formula on ergodic components and integrating \cite{AbramovInd59}. If $f$ is ergodic, then $s=1$ and $F$ is ergodic.

Abramov's time-scaling theorem \cite{Abramov59}, in the general continuous-probability-flow form proved in \cite[Section~2, Theorem]{Ito1971}, gives
\begin{equation}\label{et:time_scaling}
    h_\nu(\phi^t)=|t|h_\nu(\phi^1)\quad(t\ne0),\qquad h_\nu(\phi^0)=0.
\end{equation}
Thus for an ergodic time map $f=\phi^\tau$ with $\tau>0$,
\begin{displaymath}
    \tau h_\nu(\phi^1)=p\,h_{\nu_A}(F),\qquad
    \mathbb E_{\nu_A}N=1/p.
\end{displaymath}
For an invertible probability-preserving map and its ergodic decomposition $\nu=\int\nu_\omega\dd\lambda(\omega)$, the entropy decomposition theorem \cite[Theorem~2.8.11(b)]{PU2010} gives
\begin{equation}\label{et:ergodic_decomposition}
    h_\nu(f)=\int h_{\nu_\omega}(f)\dd\lambda(\omega).
\end{equation}
For a full-flow ergodic decomposition, the corresponding formula is
\begin{displaymath}
    h_\nu(\phi^1)=\int h_{\nu_\omega}(\phi^1)\dd\lambda(\omega).
\end{displaymath}
More generally, any measurable mixture $\eta=\int\eta_z\dd\kappa(z)$ of $f$-invariant probabilities satisfies
\begin{equation}\label{eq:entropy_decom}
    h_\eta(f)=\int h_{\eta_z}(f)\dd\kappa(z).
\end{equation}
Both extensions follow by decomposing the constituent probabilities into time-one or $f$-ergodic components, respectively, then using the uniqueness of ergodic decomposition and Tonelli's theorem.

In particular, suppose that $S$ is invariant and Borel, $a<\infty$, and every ergodic invariant probability on $S$ has entropy at most $a$. For a probability $\nu$ concentrated on $S$, almost every ergodic component is also concentrated on $S$, and the decomposition formulas give
\begin{equation}\label{eq:et:entropy_deficit}
    a-h_\nu(f)=\int\bigl(a-h_{\nu_\omega}(f)\bigr)\dd\lambda(\omega)\ge0.
\end{equation}
The same formula holds for full-flow components with $f=\phi^1$. If the component deficit is strictly positive almost everywhere, its integral is strictly positive, with no uniform gap required; moreover,
\begin{equation}\label{eq:et:entropy_equality}
    h_\nu(f)=a\quad\Longleftrightarrow\quad
    h_{\nu_\omega}(f)=a\quad\text{for almost every }\omega.
\end{equation}

Finally, the one-parameter case of Pugh--Shub's theorem \cite[Theorem~1]{PughShub71} says that an ergodic flow has at most countably many non-ergodic time maps. Its hypotheses hold here: the action on $L^2(\nu)$ is strongly continuous by continuity of the flow and density of bounded continuous functions, and $L^2(\nu)$ is separable because the phase space is a second-countable manifold. Thus, ergodic sampling times can be chosen arbitrarily small or arbitrarily large.
\section{Entropy and the critical exponent}\label{sec:entropy}
In this section, we prove that the compact topological entropy and the probability variational entropy equal the critical exponent, given that is finite:
\begin{displaymath}
    h_{\mathrm{top}}^{\,c}(\phi)=H=\delta_{\Gamma}.
\end{displaymath}
The upper bound uses the Hadamard covering structure and $\delta_{\Gamma}<\infty$; the lower bound requires only a periodic weak rank one axis.

\subsection{The upper bound for invariant probabilities}\label{sec:upper}

Let $f=\phi^1$ be the time-one map of the geodesic flow. For every $f$-invariant probability $\nu$, we prove
\begin{displaymath}
    h_\nu(f)\le\delta_{\Gamma}
\end{displaymath}
by counting orbit segments with endpoints in a fixed compact set. Distance convexity bounds the covering number for each deck transformation, and orbital growth bounds the number of transformations.

For $\widetilde v\in T^1X$, let $c_{\widetilde v}:\R\to X$ be the unit-speed complete geodesic with initial vector $\widetilde v$. Define
\begin{equation}\label{eq:up:geodesicmetric}
    d_G(\widetilde v,\widetilde w)=\int_\R d_X\left(c_{\widetilde v}(s),c_{\widetilde w}(s)\right)e^{-2|s|}\dd s.
\end{equation}
The unit-speed bound $d_X(c_{\widetilde v}(s),c_{\widetilde w}(s))\le d_X(c_{\widetilde v}(0),c_{\widetilde w}(0))+2|s|$ makes the integral finite. Symmetry and the triangle inequality follow by integration, and continuity of the integrand gives positive definiteness.

The convergence of initial vectors yields uniform convergence of geodesics on compact time intervals, and the unit-speed bound controls the tails of the integral. Conversely, if $d_G(\widetilde v_j,\widetilde v)\to0$, the $2$-Lipschitz functions $a_j(s)=d_X(c_{\widetilde v_j}(s),c_{\widetilde v}(s))$ converge uniformly to zero on every compact interval: otherwise the Lipschitz bound would give a fixed positive lower bound for their weighted integrals. The initial vectors then converge because $c'(0)=\exp_{c(0)}^{-1}(c(1))$ depends continuously on $(c(0),c(1))$. Thus $d_G$ induces the usual topology on $T^1X$.

Convexity gives a useful estimate for the standard projection $\pi_X:T^1X\to X$. If $a(s)=d_X(c_{\widetilde v}(s),c_{\widetilde w}(s))$, then $a(s)+a(-s)\ge2a(0)$, and hence
\begin{equation}\label{eq:up:basepoint_bound}
    d_G(\widetilde v,\widetilde w)=\int_0^\infty[a(s)+a(-s)]e^{-2s}\dd s\ge d_X(\pi_X\widetilde v,\pi_X\widetilde w).
\end{equation}
Since $\Gamma$ preserves $d_G$, it induces the quotient metric
\begin{displaymath}
    d(v,w)=\inf_{\gamma\in\Gamma}d_G(\widetilde v,\gamma\widetilde w),\qquad v,w\in T^1M.
\end{displaymath}
This expression is independent of the choice of the lifts and satisfies the triangle inequality because of the $\Gamma$-invariance. If $d(v,w)=0$, choose $\gamma_j$ with $d_G(\widetilde v,\gamma_j\widetilde w)\to0$. By~\eqref{eq:up:basepoint_bound}, the points $\gamma_j\pi_X\widetilde w$ lie in a fixed compact ball. Properness of the deck transformation makes $\gamma_j$ constant along a subsequence, proving $v=w$. The quotient metric induces the usual topology: convergence in a local covering chart implies convergence in $d$, while $d(v_j,v)\to0$ allows lifts with $d_G(\widetilde v_j,\widetilde v)\le d(v_j,v)+1/j$, which converge upstairs.

Moreover, \eqref{eq:up:basepoint_bound} gives $d(v,w)\ge d_M(\pi v,\pi w)$. Thus, every closed $d$-ball lies in the unit tangent bundle over a compact closed ball of $M$ and is compact, since $d$ induces the usual topology. In particular, $(T^1M,d)$ is complete and separable.

For $T\ge0$, write
\begin{displaymath}
    d_T(v,w)=\sup_{0\le t\le T}d(\phi^tv,\phi^tw).
\end{displaymath}
For any fixed pair of lifts, one has
\begin{equation}\label{eq:up:projection}
    d_T(v,w)\le\sup_{0\le t\le T}d_G(\phi^t\widetilde v,\phi^t\widetilde w).
\end{equation}

A close precursor of the compact-endpoint estimate is the argument in \cite[proof of Proposition~3, p.~40]{OtalPeigne2004}, which counts endpoint pairs in compact lifts and uses distance convexity together with a curvature lower bound in the pinched negative-curvature setting. The endpoint-net argument also uses the convexity mechanism of Manning \cite{Manning1979} and its cocompact CAT(0) formulation in \cite[Section~2]{Ricks2021}. Here, the weighted geodesic metric and the extension of segments at both ends yield the compact-endpoint estimate without a lower bound on curvature or uniform bounds on the Sasaki geometry.

\begin{lemma}\label{lem:up:endpoint}
    For compact sets $E,F\subset X$ and $\eps>0$, there is a constant $C(E,F,\eps)<\infty$ with the following property. For every $T\ge0$ and $\gamma\in\Gamma$, the unit-speed complete geodesics satisfying $c(0)\in E$ and $c(T)\in\gamma F$ are covered by at most $C(E,F,\eps)$ lifted Bowen balls of radius $\eps$ for the interval $[0,T]$. The constant is independent of $T$ and $\gamma$.
\end{lemma}

\begin{proof}
Choose $L,\eta>0$ with $2\eta+e^{-2L}<\eps$, and finite $\eta$-nets in the compact neighborhoods $N_L(E)$ and $N_L(F)$. Group the geodesics by the net balls containing $c(-L)$ and $\gamma^{-1}c(T+L)$, choosing one index when necessary. If $c,c'$ belong to the same group, synchronous distance convexity gives
    \begin{displaymath}
        d_X(c(u),c'(u))\le2\eta\qquad(-L\le u\le T+L).
    \end{displaymath}
Unit speed then gives, for $0\le t\le T$ and $s\in\R$, $d_X(c(t+s),c'(t+s))\le2\eta+2(|s|-L)_+$, so
    \begin{displaymath}
        d_G(\phi^tc,\phi^tc')
        \le2\eta+4\int_L^\infty(s-L)e^{-2s}\dd s
        =2\eta+e^{-2L}<\eps.
    \end{displaymath}
Here $\phi^tc$ denotes the initial vector of the time-shifted geodesic. One representative from each nonempty group provides the covering; the product of the two net cardinalities is independent of $T$ and $\gamma$.
\end{proof}

For later use, let $C\subset M$ be compact. Choose finitely many relatively compact evenly covered neighborhoods covering $C$, and a compact union $E\subset B(o,D)$ of lifts of their closed sub-neighborhoods. A segment of length $T$ with both base-point endpoints in $C$ lifts from $x\in E$ to $\gamma y$ with $y\in E$, so $d(o,\gamma o)\le T+2D$. Applying Lemma~\ref{lem:up:endpoint} with $E=F$ and projecting by~\eqref{eq:up:projection} therefore covers all these segments by at most
\begin{equation}\label{eq:up:returncount}
    C_\eps N_o(T+D_C),\qquad D_C=2D,
\end{equation}
$d_T$-balls of radius $\eps$, where $C_\eps=C(E,E,\eps)$ is independent of $T$. 

We combine the returning-endpoint estimate with the noncompact form of Katok's covering inequality for continuous maps of complete separable metric spaces \cite[Theorem~2.6]{Riquelme2018}.

\begin{proposition}\label{prop:up:entropy}
    Every $\phi^1$-invariant probability measure $\nu$ on $T^1M$ satisfies $h_\nu(\phi^1)\le\delta_{\Gamma}$.
\end{proposition}

\begin{proof}
First, suppose that $\nu$ is ergodic for $f=\phi^1$. Let $N_\nu(n,\eps)$ be the least number of balls for the discrete Bowen metric $\max_{0\le j<n}d(f^jv,f^jw)$ needed to cover a set of $\nu$-measure greater than $1/2$. The cited covering inequality gives
    \begin{displaymath}
        h_\nu(f)\le
        \lim_{\eps\downarrow0}\liminf_{n\to\infty}
        \frac1n\log N_\nu(n,\eps).
    \end{displaymath}
Choose a compact $C\subset M$ with $\nu(A)>3/4$, where $A=\pi^{-1}(C)$. Invariance gives $\nu(A\cap f^{-(n-1)}A)>1/2$ for every $n\ge1$. By \eqref{eq:up:returncount}, this intersection is covered by at most $C_\eps N_o(n-1+D_C)$ continuous-time Bowen balls on $[0,n-1]$, hence by the same number of discrete Bowen balls. Therefore
    \begin{displaymath}
        \liminf_{n\to\infty}\frac1n\log N_\nu(n,\eps)\le\limsup_{n\to\infty}\frac1n\log\bigl(C_\eps N_o(n-1+D_C)\bigr)\le\delta_{\Gamma}.
    \end{displaymath}
    Letting $\eps\downarrow0$ proves $h_\nu(f)\le\delta_{\Gamma}$. The entropy decomposition \eqref{et:ergodic_decomposition} gives the same conclusion for every $f$-invariant probability.
\end{proof}

\subsection{Compact invariant sets with large entropy}\label{sec:lower}

For every $b<\delta_{\Gamma}$, we construct a compact flow-invariant set whose topological entropy exceeds $b$. This compact entropy approximation is already implicit in the negative-curvature construction of \cite[\S4, pp.~37--39, especially Fact~11 on p.~38]{OtalPeigne2004}: for each finite alphabet, the separated periodic orbits remain in a common compact set, so their union has compact flow-invariant closure. Here, we adapt the construction using contracting elements and explicitly control time alignment and the multiplicity of word labels in a Bowen ball. Throughout this subsection, assume $0<\delta_{\Gamma}<\infty$, and the action is non-elementary and proper with a contracting element.

A closed subset of a geodesic metric space is called \emph{contracting} if the nearest-point projections of metric balls disjoint from it have uniformly bounded diameter. An isometry is a contracting element if its cyclic orbit is contracting and the orbit map from $\Z$ is a quasi-isometric embedding.

Fix $o\in X$. By \cite[Theorem~5.4]{BF09}, an axis of an axial isometry of a proper CAT(0) space is contracting if and only if it does not bound a flat half-plane. Thus, a periodic weak rank one geodesic supplies a contracting element: its axis and any cyclic orbit have finite Hausdorff distance, and contraction is preserved by such a change \cite[Proposition~2.4(2)]{Yang17}. A positive translation length gives the required quasi-isometric orbit map.

We apply Yang's fixed-connector construction \cite[Lemma~2.19 and the proof of Lemma~3.1]{Yang17} directly, with its ordered fellow-traveling conclusion \cite[Definition~2.6 and Proposition~2.7]{Yang17}. These statements require a proper action with a contracting element, not statistical convex co-compactness. Our non-elementary group is neither finite nor virtually cyclic, in accordance with the convention in \cite[Section~1.1]{Yang17}.

Take $\Delta=1$ and write $a(L)=\#\{\gamma:|d(o,\gamma o)-L|\le1\}$. Since $a(L)\le N_o(L+1)$, its upper exponential rate is at most $\delta_{\Gamma}$. If the rate along integer radii were smaller, a bound $a(k)\le e^{uk}$ for some $0<u<\delta_{\Gamma}$ and all large $k$, summed over the annuli covering an orbital ball, would give $N_o(r)\le C_u e^{u(r+1)}$, a contradiction. Thus there are $L_j\to\infty$ with $L_j^{-1}\log a(L_j)\to\delta_{\Gamma}$.

Let $R_0$ and the finite connector set $F_0$ be supplied by Yang's lemma. Greedy selection in the $L_j$-annulus gives an $R_0$-separated set $Z_j$ with $|Z_j|\ge a(L_j)/N_o(R_0)$: translating any selected point to $o$ bounds the number of deleted orbit points by $N_o(R_0)$. For all sufficiently large $j$, the lemma gives $B_j\subset Z_j$ with $|B_j|\ge|Z_j|/16$ and a connector $f_j\in F_0$. Taking every connector exponent equal to one in its injective extension map shows that $A_j=B_jf_j$ freely generates a positive semigroup, with
\begin{equation}\label{eq:largealphabet}
    \lim_{j\to\infty}\frac{\log|A_j|}{L_j}=\delta_{\Gamma},\qquad \ell_j:=\max_{a\in A_j}d(o,ao)\le L_j+C_0.
\end{equation}
Here $C_0=1+\max_{f\in F_0}d(o,fo)$. Indeed, $a(L_j)/(16N_o(R_0))\le |A_j|\le a(L_j)$, and the displacement bound follows from the triangle inequality.

The precise shadowing used below requires a geometric passage from Yang's admissible paths. The paths alternate between letter and connector blocks. Proposition~2.7 of \cite{Yang17} supplies an arc-length quasi-geodesic inequality and ordered points on each geodesic chord within a fixed distance $\varepsilon_0$ of the block endpoints. For each fixed alphabet $A_j$, all block lengths are bounded by some $J$. Hence, the path lies within $J+\varepsilon_0$ of its chord, while consecutive ordered chord points are at most $J+2\varepsilon_0$ apart, giving the reverse Hausdorff bound. Straightening each letter--connector pair changes this bound by a fixed amount and preserves the quasi-geodesic inequality: arc length between principal vertices decreases, and endpoints inside generator segments contribute only a bounded additive error. Thus, every finite positive-word path, and every sub-path, has uniform arc-length quasi-geodesic constants and a uniform Hausdorff bound from its chord. These constants may depend on $j$, but not on the word or its length.

We next associate complete geodesics with two-sided words in one of these alphabets.

Fix one alphabet $A=A_j$, and put
\begin{displaymath}
    \ell=\max_{a\in A}d(o,ao)>0.
\end{displaymath}
All constants in the next construction may depend on this fixed $A$. For a two-sided sequence $(a_i)_{i\in\Z}$, define group elements by $x_0=1$ and $x_{i+1}=x_i a_{i+1}$. Join successive points $x_io$ by geodesic segments, and parametrize the resulting two-sided path $q$ by arc length with $q(0)=o$.

Every finite path between vertices becomes a positive-word path after left multiplication by the inverse of its first vertex. Hence there are common constants $\lambda\ge1$ and $c\ge0$ such that
\begin{equation}\label{eq:lower_quasi-geodesic}
    \lambda^{-1}|s-t|-c\le d(q(s),q(t))\le|s-t|,
    \qquad s,t\in\R.
\end{equation}
There is also $D<\infty$ such that each finite path between vertices and its geodesic chord has a Hausdorff distance at most $D$. Allowing endpoints inside generator segments changes these constants by at most a fixed amount depending on $\ell$.

To pass to a complete geodesic, join increasingly distant negative and positive vertices by geodesic chords $c_k$. Since $q(0)$ is within $D$ of every chord, parametrize each chord so that $c_k(0)$ is the nearest point to $q(0)$. Inequality \eqref{eq:lower_quasi-geodesic} implies that the distances from $c_k(0)$ to both endpoints tend to infinity. The initial vectors $c_k'(0)$ lie over the compact ball $\overline B(o,D)$ and have a convergent subsequence. Its limit determines a complete unit-speed geodesic $c$.

For every fixed $s$, the point $q(s)$ is within $D$ of all sufficiently long chords. The corresponding nearest points on the chords stay in a bounded parameter interval, so passage to the limit gives $d(q(s),c)\le D$. Conversely, for a fixed $t$, choose $s_k$ with $d(c_k(t),q(s_k))\le D$. Then
\begin{displaymath}
    d(o,q(s_k))\le |t|+2D.
\end{displaymath}
The quasi-geodesic inequality bounds $s_k$, and a convergent subsequence gives $d(c(t),q)\le D$. We have proved
\begin{equation}\label{eq:lower_shadow}
    d_{\mathrm H}(q(\R),c(\R))\le D.
\end{equation}

We also need the orientation of this shadowing. Nearest-point projection onto $c$ is $1$-Lipschitz in the CAT(0) space $X$ \cite[Chapter~I, Corollary~5.6]{Ballmann95}. Let $r(s)$ be the signed coordinate of the projection of $q(s)$ onto $c$. Equations \eqref{eq:lower_quasi-geodesic} and \eqref{eq:lower_shadow} imply that $|r(s)|\to\infty$ as $s\to\pm\infty$. The two tails cannot approach the same end of $c$. If they did, continuity would give parameters $s_-<0<s_+$ at arbitrarily distant common projection coordinates. Then $d(q(s_-),q(s_+))\le2D$, whereas $s_+-s_-\to\infty$, contradicting \eqref{eq:lower_quasi-geodesic}. We may therefore orient $c$ so that
\begin{displaymath}
    r(s)\longrightarrow-\infty\quad(s\to-\infty),
    \qquad
    r(s)\longrightarrow+\infty\quad(s\to+\infty).
\end{displaymath}
The complete geodesic thus shadows the path with the same ordered pair of ends.

The following lemma uses the shadowing construction and a bound on the multiplicity of word labels in Bowen balls to obtain a compact flow-invariant set with the required entropy.

\begin{lemma}\label{lem:semigroup_entropy}
    Let $A\subset\Gamma$ be a finite nonempty alphabet freely generating a positive semigroup. Suppose that its two-sided positive-word paths have uniform arc-length quasi-geodesic constants and are shadowed, with their ordered ends preserved, by complete geodesics at a uniform Hausdorff distance $D$. Then there is a compact flow-invariant set $K_A\subset T^1M$ such that
    \begin{displaymath}
        h_{\mathrm{top}}(\phi^1|K_A)\ge\frac{\log|A|}{\ell},
        \qquad \ell=\max_{a\in A}d(o,ao).
    \end{displaymath}
\end{lemma}

\begin{proof}
    Let $E\subset M$ be the union of the projections of the finitely many generator segments $[o,ao]$, $a\in A$. Every segment in every positive-word path projects into $E$. The projection of each shadowing geodesic is consequently contained in
    \begin{displaymath}
        E_D=\{x\in M:d(x,E)\le D\}.
    \end{displaymath}
    This set is compact because $E$ is compact. Let $K_A$ be the closure of all unit tangent vectors on the projected shadowing geodesics. It lies in the compact unit tangent bundle over $E_D$ and is flow-invariant by continuity of the flow. Fix a letter $a_*\in A$. Extend each positive word $w=a_1\cdots a_n$ by repeating $a_*$ indefinitely on both sides. This gives a two-sided path $q_w$ with zeroth vertex $o$. Its vertex $wo$ occurs at the arc length parameter
    \begin{displaymath}
        0<s_w=\sum_{i=1}^n d(o,a_io)\le n\ell.
    \end{displaymath}
    Choose a complete shadowing geodesic $c_w$ with the corresponding orientation, initially with any time origin. Denote by $r_w(s)$ the signed coordinate of the nearest-point projection of $q_w(s)$ onto $c_w$. If $s>0$ and $r_w(s)<r_w(0)$, the limit $r_w(t)\to+\infty$ gives an $s'>s$ with $r_w(s')=r_w(0)$. The two corresponding path points have the same projection and hence satisfy
    \begin{displaymath}
        d(q_w(0),q_w(s'))\le2D,
        \qquad s'\le B:=\lambda(2D+c).
    \end{displaymath}
    As projection is $1$-Lipschitz, $r_w(0)-r_w(s)\le s\le B$. Thus $r_w(s)\ge r_w(0)-B$ for every $s\ge0$. Set the time-zero point of $c_w$ at the old coordinate $r_w(0)-B$. This choice yields
    \begin{displaymath}
        d(o,c_w(0))\le D':=D+B.
    \end{displaymath}
    The endpoint coordinate $t_w=r_w(s_w)-r_w(0)+B$ then satisfies
    \begin{equation}\label{eq:word_time}
        0\le t_w\le n\ell+B,
        \qquad d(wo,c_w(t_w))\le D.
    \end{equation}
    Write $T_n=n\ell+B$. Every endpoint can now be read in the same positive time interval $[0,T_n]$, and every chosen lift starts in the same compact ball $\overline B(o,D')$.

    To count how many words can lie in one Bowen ball, use the Sasaki metric $d_{\rm Sas}$ on $T^1M$. Its standard projection is $1$-Lipschitz. On $K_A$ define
    \begin{displaymath}
        d_T^{\rm Sas}(v,w)=\max_{0\le t\le T}
        d_{\rm Sas}(\phi^tv,\phi^tw).
    \end{displaymath}
    The compactness of $E_D$ gives $i_0=\inf_{x\in E_D}\operatorname{inj}_M(x)>0$. Fix $0<\eps<i_0/4$.

    Choose a reference word $w_0$ and its selected lift $c_0=c_{w_0}$. Suppose that the initial vector of the projection of $c_w$ lies in the $d_{T_n}^{\rm Sas}$-ball of radius $\eps$ centered at the reference vector. At every time $t\in[0,T_n]$, there is a deck transformation $\gamma(t)$ satisfying
    \begin{displaymath}
        d_X(c_0(t),\gamma(t)c_w(t))<\eps.
    \end{displaymath}
    It is unique. If two elements $\gamma_1,\gamma_2$ satisfied this inequality, then
    \begin{displaymath}
        d_X(c_w(t),\gamma_1^{-1}\gamma_2c_w(t))<2\eps<2i_0,
    \end{displaymath}
    which contradicts the lower bound of twice the injectivity radius for the displacement of a nontrivial deck transformation. Continuity and the strict inequality make $\gamma(t)$ locally constant. It is therefore a single element $\gamma$ on the connected interval $[0,T_n]$.

    The estimate at time zero and the two initial basepoints bounds give
    \begin{displaymath}
        d(o,\gamma o)\le2D'+\eps.
    \end{displaymath}
    Thus, all possible alignment elements lie in the finite set
    \begin{displaymath}
        F=\{\gamma\in\Gamma:d(o,\gamma o)\le2D'+\eps\}.
    \end{displaymath}
    For each such fixed $\gamma$, equation \eqref{eq:word_time} puts $\gamma wo$ in the $(D+\eps)$-neighborhood of the reference segment $c_0([0,T_n])$.

    Any ball $B(x,R)$ contains at most $N_o(2R)$ orbit points. If it contains $\beta o$, left multiplication by $\beta^{-1}$ maps all its orbit points into $B(o,2R)$. The action of the deck transformation is free, so distinct group elements give distinct orbit points. Cover the reference segment by at most $T_n+2$ balls of radius $1$. Its $(D+\eps)$-neighborhood contains at most
    \begin{displaymath}
        (T_n+2)N_o(2D+2\eps+2)
    \end{displaymath}
    orbit points. For each fixed $\gamma$, the map $w\mapsto\gamma wo$ is injective because the positive semigroup is free. If $\bar c_w$ denotes the projected geodesic, we have proved
    \begin{equation}\label{eq:polynomial_fiber}
        \#\{w:\ d_{T_n}^{\rm Sas}(\bar c_w'(0),\bar c_0'(0))<\eps\}
        \le |F|(T_n+2)N_o(2D+2\eps+2).
    \end{equation}
    This counts word labels even when their projected initial vectors coincide.

    Apply a greedy selection to the $|A|^n$ labeled initial vectors. Choose one vector, delete every word whose vector lies in its Bowen $\eps$-ball, and continue. Each selection deletes at most the number in \eqref{eq:polynomial_fiber}. The selected vectors are $(T_n,\eps)$-separated and satisfy
    \begin{displaymath}
        s(T_n,\eps,K_A)\ge
        \frac{|A|^n}{|F|(T_n+2)N_o(2D+2\eps+2)}.
    \end{displaymath}
    Taking logarithms, dividing by $T_n$, and letting $n\to\infty$ gives
    \begin{displaymath}
        \limsup_{T\to\infty}\frac1T\log s(T,\eps,K_A)
        \ge\frac{\log|A|}{\ell}.
    \end{displaymath}
    The entropy formula in Section~\ref{sec:preliminaries} gives the claim.
\end{proof}

For each fixed alphabet $A_j$, all constants in the preceding proof are independent of the word length. After taking the entropy limit, let $j\to\infty$ in \eqref{eq:largealphabet} to obtain
\begin{equation}\label{eq:lower_compact_entropy}
    h_{\mathrm{top}}^{\,c}(\phi)=\sup_{\substack{K\subset T^1M~\text{compact}\\\phi^tK=K\ (t\in\R)}}h_{\mathrm{top}}(\phi^1|K)\ge\limsup_{j\to\infty}\frac{\log|A_j|}{L_j+C_0}=\delta_{\Gamma}.
\end{equation}

Time averaging preserves entropy by \eqref{eq:entropy_decom}. If $K$ is compact and flow invariant and $\eta$ is $\phi^1|_K$-invariant, then $\eta_t=(\phi^t)_*\eta$ is a $1$-periodic family of $\phi^1$-invariant probabilities. Its average is flow invariant, supported on $K$, and satisfies
\begin{equation}\label{eq:lower:time_average}
    \begin{aligned}
        \nu&=\int_0^1(\phi^t)_*\eta\dd t,\\
        h_\nu(\phi^1)&=\int_0^1h_{\eta_t}(\phi^1)\dd t=h_\eta(\phi^1).
    \end{aligned}
\end{equation}
The invariance follows by translating an integral of a $1$-periodic family, and the last equality follows because $\phi^t$ is a measure-theoretic isomorphism from $(K,\eta,\phi^1)$ to $(K,\eta_t,\phi^1)$. 

\begin{theorem}\label{thm:H_delta}
    The supremum of the topological entropies on compact flow-invariant sets and the supremum of the entropies of invariant probabilities satisfy
    \begin{displaymath}
        h_{\mathrm{top}}^{\,c}(\phi)=H=\delta_{\Gamma}.
    \end{displaymath}
    More precisely, for every $b<\delta_{\Gamma}$ there is a compactly supported flow-invariant probability $\nu$ with $h_\nu(\phi^1)>b$.
\end{theorem}

\begin{proof}
    Section~\ref{sec:preliminaries}, Proposition~\ref{prop:up:entropy}, and \eqref{eq:lower_compact_entropy} give
    \begin{displaymath}
        \delta_{\Gamma}\le h_{\mathrm{top}}^{\,c}(\phi)\le H\le\delta_{\Gamma}.
    \end{displaymath}
    For $b<\delta_{\Gamma}$, choose a compact flow-invariant set $K$ with $h_{\mathrm{top}}(\phi^1|K)>b$. The compact variational principle provides a $\phi^1|K$-invariant probability $\eta$ with $h_\eta(\phi^1)>b$, and~\eqref{eq:lower:time_average} gives a flow-invariant probability supported on $K$ with the same entropy.
\end{proof}

\section{Maximal entropy measures and the finiteness criterion}\label{sec:uniqueness}

The same unstable partitions yield the entropy of the finite BMS measure and establish rigidity for regular invariant probabilities with entropy $\delta_{\Gamma}$. We then exclude singular maximizing probabilities and complete the proof of Theorem~\ref{thm:main}.

\subsection{Unstable partitions and reference measures}\label{sec:unstable-tools}

We adapt the partition and conditional-entropy methods of Otal and Peign\'e \cite[Propositions~1 and~4, Lemma~8, and \S3]{OtalPeigne2004}; see also \cite[Proposition~6.4 and Lemma~6.6]{PPS2015}. Recurrent compact windows and finite-entropy return and deck coding, together with a flow-coordinate refinement, provide the required partitions.

Fix an arbitrary flow-ergodic invariant probability $\nu$ with $\nu(\Reg)=1$. Write $Y=T^1M$ and use $\pi$ for the standard projections on the cover and quotient. By \cite[Theorem~1]{PughShub71}, choose $\tau>0$ so that $f=\phi^\tau$ is ergodic. Use the Borel-partition and entropy conventions of Section~\ref{sec:entropy_tools}, writing $\bar\nu_\alpha$ for the atom distribution, $\nu_A$ for the conditional probabilities, and $\alpha^-=\bigvee_{j\geq0}f^j\alpha$ for the present and entire past. 

Fix $o\in X$ and use the Busemann convention and normalized conformal density of Section~\ref{sec:geometry}. Write $v^-,v^+$ for the endpoints of $v\in T^1X$. Horospherical strong unstable and strong stable horospheres are determined by $(v^-,\beta_{v^-}(\pi v,o))$ and $(v^+,\beta_{v^+}(o,\pi v))$, respectively. A plaque is a complete local leaf slice in a foliation chart. Leaf-path lengths are horizontal (base-projection) lengths unless stated otherwise; Sasaki comparisons will be made on specified compact sets.

By \cite[Chapter~IV, Propositions~2.8 and~3.2]{Ballmann95}, the Green operators are continuous, Busemann functions are $C^2$, and the horospheres' unit normal graphs form continuously varying $C^1$ leaves of dimension $n-1$, tangent to $(Z,U^sZ)$ and $(Z,U^uZ)$. Convexity gives $U^s\leq0\leq U^u$; time reversal exchanges the two cases. In a parallel normal frame, the unstable Jacobi tensor satisfies
\begin{displaymath}
    A(0)=I,\qquad A'=U^u(\phi^tv)A,\qquad A''+R_tA=0.
\end{displaymath}
It is invertible for every $t$. Since $U^u\geq0$, the Jacobi field $J(t)=A(t)z$ satisfies
\begin{displaymath}
    \frac{d}{dt}|J(t)|^2
    =2\langle U^u(\phi^tv)J(t),J(t)\rangle\geq0.
\end{displaymath}
Consequently, the backward flow does not increase the horizontal length of any unstable path.

We use the parallel-Jacobi-field criterion of \cite[Chapter~IV, Lemma~2.6]{Ballmann95} to extract a positive window along each regular orbit.

\begin{lemma}\label{lem:rigidity:positive_window}
    For every $v\in\Reg$, there is a positive integer $k$ for which $U^u(\phi^kv)$ is positive definite. If $\nu(\Reg)=1$, then the open set $\{U^u>0\}$ has positive $\nu$ measure.
\end{lemma}
\begin{proof}
    Identify the normal spaces along $c_v$ by parallel transport and set
    \begin{displaymath}
        Q(t)=A(t)^*U^u(\phi^tv)A(t)=A(t)^*A'(t).
    \end{displaymath}
    This is a positive semidefinite self-adjoint matrix. The Jacobi equation gives
    \begin{displaymath}
        Q'(t)=A'(t)^*A'(t)-A(t)^*R_tA(t)\geq0,
    \end{displaymath}
    because $R_t$ is negative semidefinite. Suppose that $U^u(\phi^kv)$ is singular for every positive integer $k$. Then $\ker Q(k)$ is a decreasing sequence of nonzero subspaces admitting a nonzero vector $z$ in the intersection. For any $t\in\R$, choose $k\geq\max\{1,t\}$; the inequalities $0\leq Q(t)\leq Q(k)$ imply $\langle Q(t)z,z\rangle=0$. For $J(t)=A(t)z$ we have
\begin{displaymath}
    \frac{d}{dt}|J(t)|^2=2\langle Q(t)z,z\rangle=0.
\end{displaymath}
Thus $J$ is a nonzero normal Jacobi field of constant norm, hence is parallel by \cite[Chapter~IV, Lemma~2.6]{Ballmann95}, contrary to regularity.

    Finally,
    \begin{displaymath}
        \Reg\subset\bigcup_{k\geq1}\phi^{-k}\{U^u>0\}.
    \end{displaymath}
    Flow invariance and $\nu(\Reg)=1$ imply that $\nu\{U^u>0\}>0$.
\end{proof}

Let $\widetilde\Reg\subset T^1X$ be the lift of the regular set and put $t(v)=\beta_{v^-}(\pi v,o)$. Then
\begin{equation}\label{eq:rigidity-hopfmap}
    \Phi(v)=(v^-,v^+,t(v))
\end{equation}
is a homeomorphism from $\widetilde\Reg$ onto the open set $\mathcal E_{\Reg}\times\R$, where $\mathcal E_{\Reg}$ consists of the ordered endpoint pairs of regular geodesics. We use Hopf coordinates on evenly covered neighborhoods in $\Reg$. A Hopf box is such a neighborhood whose chosen lift corresponds to a product $I_-\times I_+\times J$ in these coordinates. Indeed, restrict the zero-width Hopf homeomorphism of \cite[Section~5]{Link2018} to rank one geodesics and negate its last coordinate. The endpoint neighborhoods in \cite[Section~3]{Link2018} make $\mathcal E_{\Reg}$ open, while flow invariance gives the product with $\R$. Fixing $(v^-,t)$ therefore makes the positive endpoint a topological coordinate on local strong unstable plaques.

Choose $v_0\in\supp\nu\cap\Reg\cap\{U^u>0\}$ by Lemma~\ref{lem:rigidity:positive_window}. Choose $0<\eta<\tau/3$. By continuity, after decreasing $\eta$ and the neighborhood of $v_0$, there are a compact neighborhood $C$, a relatively compact open neighborhood $C_1\supset C$, and $a>0$ such that
\begin{equation}\label{eq:rigidity-negwindow}
    U^u(\phi^sw)\geq aI
    \qquad(w\in C_1,\ -\eta\leq s\leq0).
\end{equation}

Let $K_1$ be the closed unit neighborhood of $\pi C$. Completeness makes $K_1$ compact, so
\begin{displaymath}
    B_1=\sup_{T^1K_1}\|U^u\|<\infty.
\end{displaymath}
Set $b=\dist_{\mathrm{Sas}}(C,C_1^c)>0$ and choose
\begin{displaymath}
    0<a_0<\min\{1,b/\sqrt{1+B_1^2}\}.
\end{displaymath}
An unstable path starting in $C$ with horizontal length less than $a_0$ has base points in $K_1$ and Sasaki length at most $\sqrt{1+B_1^2}a_0<b$, so it stays in $C_1$. Throughout its backward image over $[-\eta,0]$, \eqref{eq:rigidity-negwindow} gives $\frac{d}{dt}|J|\geq a|J|$. This window multiplies horizontal length by at most
\begin{displaymath}
    \theta=e^{-a\eta}<1.
\end{displaymath}
Outside these windows, backward horizontal length remains non-increasing.

Choose a relatively compact, evenly covered neighborhood $W\subset\Reg$ equipped with Hopf coordinates with closure in $C$, inside a larger continuous $C^1$ foliation chart. Shrinking $W$ and using coordinate bounds on compact parameter sets, arrange that any two points of a plaque are joined in the larger leaf disk by a path of horizontal length less than $a_0$. Use this $C^1$ leaf chart for differentiable estimates. Choose $r_0>0$ and an open set $C_0$ such that
\begin{align*}
&\overline{B_{\mathrm{Sas}}(v_0,2r_0)}\subset W,\\
&v_0\in C_0,\qquad
\overline{C_0}\subset B_{\mathrm{Sas}}(v_0,r_0),
\qquad \nu(C_0)>0.
\end{align*}
All these choices are fixed while the radius of the return ball varies in $(r_0,2r_0)$.

Let $N_-(v,k)$ count the visits of the center to $C_0$ in the past $k$ sampling times. Backward nonexpansion and the disjoint contraction windows give, simultaneously for every unstable path $I$ starting at $v$ with horizontal length less than $a_0$,
\begin{displaymath}
 \operatorname{length}_h(\phi^{-k\tau}I)
 \leq\theta^{N_-(v,k)-O(1)}\operatorname{length}_h(I).
\end{displaymath}
By Birkhoff ergodic theorem, $N_-(v,k)/k\to\nu(C_0)>0$ almost everywhere. Fixing $0<\lambda<-\nu(C_0)\log\theta$ and absorbing the finitely many initial deviations gives a finite $A(v)$ such that
\begin{displaymath}
 \operatorname{length}_h(\phi^{-k\tau}I)
 \leq A(v)e^{-\lambda k}\operatorname{length}_h(I),\qquad k\geq0.
\end{displaymath}
The same count applies to a path of length less than $a_0$ starting at a future iterate and followed backward to $v$, with only a bounded number of endpoint windows omitted. All estimates hold simultaneously for the whole indicated plaque, not merely for typical points in it.

Let
\begin{displaymath}
    V=B_{\mathrm{Sas}}(v_0,r),\qquad r\in(r_0,2r_0).
\end{displaymath}
The estimates below hold for every such radius, to be fixed by a boundary estimate. Fix an inverse covering chart $L:W\to\widetilde W$, put $\widetilde V=L(V)$, and choose $D<\infty$ with $\pi\widetilde V\subset B(o,D)$ for all allowed radii.

Let $B$ be the $f$-invariant Borel set of points that visit $V$ infinitely often in both time directions. Ergodicity and $\nu(V)>0$ give $\nu(B)=1$. For $v\in V\cap B$, define
\begin{displaymath}
    n(v)=\min\{n\geq1:f^{-n}v\in V\}
\end{displaymath}
and the unique deck increment $\gamma(v)\in\Gamma$ by
\begin{equation}\label{eq:rigidity-deck}
    \phi^{-n(v)\tau}L(v)=\gamma(v)L(f^{-n(v)}v).
\end{equation}
Let $\cP$ record $(n(v),\gamma(v))$ on $V\cap B$, with one complementary atom containing $Y\setminus B$. The Borel return time, unique deck increment, Borel relations~\eqref{eq:rigidity-deck}, and countability of $\Gamma$ make $\cP$ a countable Borel partition. All geometric atoms below remain defined on the full space.

To prove finite entropy, note that \eqref{eq:rigidity-deck} implies
\begin{displaymath}
    d(o,\gamma(v)o)\leq n(v)\tau+2D.
\end{displaymath}
For $s>\delta_{\Gamma}$, the orbital bound~\eqref{eq:fb:orbitalcount} bounds the number $M_n$ of possible increments at return time $n$ by
\begin{displaymath}
    M_n\leq C_{s,D}e^{sn\tau}\qquad(C_{s,D}\geq1).
\end{displaymath}

Put $p=\nu(V)$. Kac's formula~\eqref{et:kac} for $f^{-1}$ gives $\mathbb E_{\nu|_V/p}n=p^{-1}$. The conditional label estimate~\eqref{et:finite_labels} and the cost of recording membership in $V$ give
\begin{displaymath}
    H_\nu(\cP)
    \leq H(p,1-p)+p\log C_{s,D}+\log2+s\tau<\infty.
\end{displaymath}

Put $\cP^-:=\bigvee_{k\geq0}f^k\cP$. If $v\in B$ and $w\in\cP^-(v)$, a common past visit to $V\cap B$ gives $w\in B$. Starting from the most recent visit, the equal-return labels in \eqref{eq:rigidity-deck} give coherent lifts with the same accumulated deck word at every earlier return. Their base-point distance is at most $\diam(\pi\widetilde V)$ at returns and, by convexity, between returns; unit speed bounds the remaining finite interval. Therefore, the lifts have the same negative endpoint. Undoing the common deck word gives equal canonical negative endpoints at every return.

We adapt the unstable-partition construction of \cite[Proposition~1, pp.~24--26]{OtalPeigne2004}; see also \cite[Proposition~6.4 and its following remark]{PPS2015}. Recurrent-window contraction replaces uniform negative-curvature contraction and applies to the complete Borel atoms.

For $v\in V$, let $\widehat\zeta(v)$ be the intersection with $V$ of its complete plaque in $W$, determined by its negative endpoint and negative Busemann time. All intersection components belong to one atom and retain the uniformly short connecting paths in the larger leaf disk. Use one complementary atom on $Y\setminus V$, and define
\begin{displaymath}
    \zeta=\bigvee_{k\geq0}f^k\widehat\zeta,
    \qquad \zeta_n=f^{-n}\zeta\quad(n\geq0).
\end{displaymath}
On all of $Y$, $\widehat\zeta$ and $\zeta$ are respectively the fiber partitions of the Borel maps
\begin{displaymath}
    \ell(v)=
    \begin{cases}
        \left(L(v)^-,\beta_{L(v)^-}(\pi L(v),o)\right),&v\in V,\\
        *,&v\notin V,
    \end{cases}
    \qquad
    c(v)=\left(\ell(f^{-k}v)\right)_{k\geq0}.
\end{displaymath}
The target $\mathcal Z$ of $\ell$ is the local $(\xi,t)$ parameter space with one additional symbol $*$; it and $\mathcal Z^{\{0,1,\ldots\}}$ are standard Borel. In particular, $w\in\zeta(v)$ is the Borel relation $c(w)=c(v)$.

The partitions $\zeta_n$ successively refine one another. For a point with a past visit to $V$, put
\begin{displaymath}
    N(v)=\min\{k\geq0:f^{-k}v\in V\}.
\end{displaymath}
This is Borel and constant on each such $\zeta$ atom. Moreover,
\begin{equation}\label{eq:rigidity-plaque-container}
    \zeta(v)\subset f^{N(v)}\widehat\zeta(f^{-N(v)}v).
\end{equation}
The right-hand side lies within a finite flow image of a relatively compact embedded leaf disk. Thus, almost every atom is relatively compact in its leaf.

\begin{lemma}\label{lem:rigidity-subordinate}
    The radius $r$ can be chosen so that $\zeta$ is subordinate to strong unstable leaves for $\nu$: for almost every $v$, the atom $\zeta(v)$ contains a leafwise open neighborhood of $v$ and is relatively compact in its leaf. Furthermore,
    \begin{equation}\label{eq:rigidity-actual-diameter}
        \diam_{\mathrm{Sas}}\left(\zeta_n(v)\right)\longrightarrow0
        \quad\text{for }\nu\text{-almost every }v,
    \end{equation}
    where the diameter is taken over the complete Borel atom.
\end{lemma}
\begin{proof}
    If an unstable path of horizontal length at most $a_0$ meets $\partial V$, all its base points lie in the closed $a_0$ neighborhood of $\pi\overline V$. Hopf--Rinow makes this neighborhood compact. The Green-operator continuity from \cite[Chapter~IV, Proposition~2.8]{Ballmann95} bounds $\|U^u\|$ on its unit tangent bundle, so the tangent representation $(Z,U^uZ)$ bounds Sasaki length by $C$ times horizontal length. This also applies when the center is outside $V$, and $C$ can be chosen uniformly for $r\in(r_0,2r_0)$.

    For every $\eps>0$, Tonelli's theorem gives
    \begin{displaymath}
        \int_{r_0}^{2r_0}\sum_{n\geq0}
        \nu\{|d_{\mathrm{Sas}}(v,v_0)-r|\leq e^{-\eps n}\}\dd r
        \leq2\sum_{n\geq0}e^{-\eps n}<\infty.
    \end{displaymath}
    Choose $\eps=\lambda/2$ and a radius for which this series is finite and $\nu(\partial V)=0$. Invariance and Borel--Cantelli imply that, for almost every $v$ and all sufficiently large $n$,
    \begin{displaymath}
        \dist_{\mathrm{Sas}}(f^{-n}v,\partial V)>e^{-\eps n}.
    \end{displaymath}
    Here $\dist(x,\partial V)\geq|d(x,v_0)-r|$ follows from the $1$-Lipschitz property of the distance function.

    For a typical $v$, any unstable path from $v$ of horizontal length less than a sufficiently small $b(v)>0$ has length at most $A(v)b(v)e^{-\lambda n}$ at its $n$th past iterate. A crossing of $\partial V$ at that time would put the center within $CA(v)b(v)e^{-\lambda n}$ of the boundary, contradicting the preceding lower bound for all large $n$. Shrink $b(v)$ to handle the finitely many earlier times; the typical orbit never meets $\partial V$.

    At every past time, the path and its center therefore have the same membership in $V$. When they lie in $V\subset W$, the path lifts continuously under the same covering chart and has a constant negative endpoint and a negative Busemann height. It belongs to the same canonical plaque. Otherwise, it stays in the complementary atom. The union of these sufficiently short unstable paths is a leafwise open neighborhood of $v$ contained in $\zeta(v)$.

    To prove \eqref{eq:rigidity-actual-diameter}, take future return times $n_j\to\infty$ for the center. Then
    \begin{displaymath}
        \zeta_{n_j}(v)\subset f^{-n_j}\widehat\zeta(f^{n_j}v).
    \end{displaymath}
    Every point of the plaque on the right before applying $f^{-n_j}$ can be joined to its center by a path of horizontal length less than $a_0$. Pull these entire paths backward to time zero along the flows. The center visits $C_0$ with positive frequency during the segment, and each visit contributes a factor $\theta$. Therefore, their horizontal lengths tend uniformly to zero. In a compact neighborhood of the fixed center $v$, the local Sasaki comparison then makes their Sasaki lengths tend to zero. This proves the assertion along $n_j$ for the entire containing plaque. Since the atoms $\zeta_n(v)$ are nested, the assertion follows for every $n$.
\end{proof}

Intersect the full-measure Borel sets just obtained with all their integer translates. We may then use subordination, both visit frequencies, endpoint recurrence, and complete-atom contraction simultaneously at every integer time. Only the centers are restricted: the geometric atoms remain defined on the whole phase space.

A finite additional partition removes the remaining ambiguity in the flow direction.

Choose a small rectangular Hopf tube $U$ with closure in $V$, whose time coordinate satisfies $-a_*<\vartheta<a_*$ and which contains $\{\phi^sv_0:|s|<a_*\}$. The coordinate $\vartheta$ differs from $t$ by a constant and increases at unit speed along the flow. Set
\begin{displaymath}
    \mathcal A=\{U,U^c\},\qquad \cQ=\cP\vee\mathcal A.
\end{displaymath}
Then $H_\nu(\cQ)\leq H_\nu(\cP)+\log2<\infty$.

Put $\cQ^-:=\bigvee_{k\geq0}f^k\cQ$. Moreover, $\cQ^-$ refines $\zeta$ modulo $\nu$. To see this, take a typical $v$ and any $w\in\cQ^-(v)$. The coherent-lift observation gives a common negative endpoint, while the Busemann-time difference $b=t(w)-t(v)$ is unchanged by the flow and common deck transformations. For every $m$ with $1/m<2a_*$, the open end strips
\begin{displaymath}
 U_m^+=\{x\in U:\vartheta(x)>a_*-1/m\},\qquad
 U_m^-=\{x\in U:\vartheta(x)<-a_*+1/m\}
\end{displaymath}
have positive measure: each contains a point of the orbit segment of $v_0\in\supp\nu$, and the support is flow invariant. A typical $v$ visits both strips arbitrarily far in the past. Equality of the $\mathcal A$-names places $w$ in $U$ at those times; the canonical time difference is still $b$, so $-1/m<b<1/m$. Hence $b=0$. At every past visit to $V$, the two points therefore lie in the same strong unstable plaque, and otherwise in the complementary atom. Thus $w\in\zeta(v)$.

The finite-entropy past and complete-atom contraction give the following entropy identity, corresponding to \cite[Proposition~4, pp.~27--31]{OtalPeigne2004} and \cite[Proposition~6.4]{PPS2015}.

\begin{proposition}\label{prop:rigidity-fullentropy}
    For the partition $\zeta$ constructed above,
    \begin{displaymath}
        h_\nu(f)=H_\nu(f^{-1}\zeta\mid\zeta),
        \qquad
        H_\nu(\zeta_n\mid\zeta)=nh_\nu(f).
    \end{displaymath}
\end{proposition}
\begin{proof}
Put $c=H_\nu(f^{-1}\zeta\mid\zeta)$. Choose increasing finite partitions $E_\ell$ generating $\zeta$. Put $E_\ell^-:=\bigvee_{k\geq0}f^kE_\ell$. Since $f\zeta$ is coarser than $\zeta$, so is $E_\ell^-$, and
\begin{displaymath}
 H_\nu(f^{-1}E_\ell\mid\zeta)
 \leq H_\nu(f^{-1}E_\ell\mid E_\ell^-)
 =h_\nu(f,E_\ell)\leq h_\nu(f).
\end{displaymath}
Monotone convergence gives $c\leq h_\nu(f)\leq\tau\delta_{\Gamma}<\infty$, where the final bound follows from Proposition~\ref{prop:up:entropy} and time scaling. The nested chain rule and invariance give $H_\nu(\zeta_n\mid\zeta)=nc$.

For any finite partition $\mathcal F$, put $\mathcal R=\mathcal F\vee\cQ$. Then $H_\nu(\mathcal R)<\infty$, and $\mathcal R^-$ refines $\zeta$. The past-entropy identity, conditioning, and the chain rule yield
\begin{displaymath}
 nh_\nu(f,\mathcal R)
 \leq nc+\sum_{j=0}^{n-1}H_\nu(\mathcal R\mid\zeta_j).
\end{displaymath}
Complete-atom contraction makes the conditional expectations of a countable determining family of bounded Lipschitz functions converge to the functions themselves. Thus the increasing partitions $\zeta_j$ generate the point partition modulo $\nu$, and \eqref{et:conditioning_limit} gives $H_\nu(\mathcal R\mid\zeta_j)\to0$. Divide the displayed inequality by $n$, take the Ces\`aro limit, and then the supremum over $\mathcal F$ to obtain $h_\nu(f)\leq c$. This proves both identities.
\end{proof}

As in \cite[\S3, p.~34]{OtalPeigne2004} and \cite[Lemma~6.6]{PPS2015}, we compare conditional probabilities with normalized Patterson--Sullivan measures on the candidate probability's subordinate atoms, without assuming finiteness of $m$. We verify the kernel on the present complete Borel atoms.

Following the Patterson--Sullivan leaf construction of \cite[\S1, pp.~21--22]{OtalPeigne2004}, also given as the zero-potential case in \cite[\S3.9, equations~(48)--(51)]{PPS2015}, we use the positive endpoint as a local coordinate on the regular part of a geometric unstable horospherical leaf $H$ in $T^1X$ and define
\begin{equation}\label{eq:rigidity-leafmeasure}
    \dd m_H^u(w)
    =e^{\delta_{\Gamma}\beta_{w^+}(o,\pi w)}\dd\sigma_o(w^+).
\end{equation}
The endpoint map is open and injective on regular plaques, and the coordinates and weight agree on overlapping leaf charts, yielding consistent Borel leaf measures. For deck equivariance, put $p=\pi w$ and $\eta=w^+$; then
\begin{align*}
    \beta_{\gamma\eta}(o,\gamma p)
&=\beta_\eta(\gamma^{-1}o,p)
=\beta_\eta(\gamma^{-1}o,o)+\beta_\eta(o,p),\\
\frac{\dd\sigma_{\gamma^{-1}o}}{\dd\sigma_o}(\eta)
&=e^{-\delta_{\Gamma}\beta_\eta(\gamma^{-1}o,o)}.
\end{align*}
The two additional factors cancel in the pulled-back endpoint integral, giving $m^u_{\gamma H}(\gamma E)=m_H^u(E)$ and consistent measures on quotient leaf charts.

The positive endpoint is unchanged by the flow, while $\beta_{w^+}(o,\pi \phi^tw)$ increases by $t$. Consequently,
\begin{equation}\label{eq:rigidity-leafscale}
    m_{\phi^tH}^u(\phi^tE)=e^{\delta_{\Gamma} t}m_H^u(E).
\end{equation}
Equivalently, $(\phi^t)_*m_H^u=e^{-\delta_{\Gamma} t}m_{\phi^tH}^u$.

A typical $v$ returns to a compact set in both directions; its lifted return points stay a bounded distance from orbit points $\gamma_j o$, so both endpoints are radial limit points in $\Lambda_\Gamma$. Since $\supp\sigma_o=\Lambda_\Gamma$, every open unstable neighborhood of $v$ has positive reference measure. Boundedness of the Busemann weight on relatively compact lifted plaques, subordination, and \eqref{eq:rigidity-plaque-container} give
\begin{equation}\label{eq:rigidity-rho}
    \rho(v):=m^u_{W^u(v)}(\zeta(v))\in(0,\infty)
    \quad\text{for }\nu\text{-almost every }v.
\end{equation}
As in \cite[Fact~9, pp.~34--35]{OtalPeigne2004}, the complete Borel code will show that $\rho$ and the normalized restrictions
\begin{displaymath}
    m_A^0=\frac{m^u|_A}{m^u(A)},\qquad A=\zeta(v),
\end{displaymath}
are measurable, with restrictions taken on the complete Borel atoms in $Y$.

Let $\mathcal C=\mathcal Z^{\{0,1,\ldots\}}$ be the standard Borel coding space of $c$. For $z\in\mathcal C$ having an entry different from $*$, let $k(z)$ be the first such index and write $z_{k(z)}=(\xi_z,t_z)$. These are Borel functions of $z$. Let $\Psi$ be the inverse of the Hopf coordinate map on $\widetilde W$ and $\operatorname{pr}:T^1X\to T^1M$ the covering projection. Wherever $\Psi(\xi_z,\eta,t_z)$ is defined, put
\begin{displaymath}
    w_z(\eta)=f^{k(z)}\operatorname{pr}\Psi(\xi_z,\eta,t_z).
\end{displaymath}
For a Borel set $E\subset Y$ and $k(z)=k$, define
\begin{align}
    K(z,E)=e^{\delta_{\Gamma} k\tau}\int
 &\indic_{\widetilde V}\left(\Psi(\xi_z,\eta,t_z)\right)
 \indic_{\{c(w_z(\eta))=z\}}
 \indic_E\left(w_z(\eta)\right)\notag\\[-2mm]
 &\qquad\cdot
 e^{\delta_{\Gamma}\beta_\eta(o,\pi\Psi(\xi_z,\eta,t_z))}
 \dd\sigma_o(\eta).
 \label{eq:rigidity-kernel-explicit}
\end{align}
Set the integrand to zero outside the domain of $\Psi$, and set $K(z,E)=0$ for the code whose entries are all $*$.

Joint Borel measurability of the integrand makes $K(z,E)$ Borel in $z$. When $K(z,\cdot)$ is nonzero, the bounded weight on $\widetilde V$ and finite $k(z)$ make it finite; the equality-of-codes indicator supports it on $c^{-1}(z)$. If $z=c(v)$ and $N(v)<\infty$, then $k(z)=N(v)$, so \eqref{eq:rigidity-plaque-container} and \eqref{eq:rigidity-leafscale} give
\begin{displaymath}
    K(c(v),E)=m^u_{W^u(v)}\left(E\cap\zeta(v)\right),
    \qquad \rho(v)=K(c(v),Y).
\end{displaymath}
Thus $\rho$ is measurable. Normalize $K(z,\cdot)$ on $\{z:0<K(z,Y)<\infty\}$ and assign any fixed probability on the complement. This gives the claimed measurable probability kernel, with the normalization valid almost everywhere by \eqref{eq:rigidity-rho}. Deck equivariance and \eqref{eq:rigidity-leafscale} ensure the independence of the chosen past visit.

The reference-mass calculation follows \cite[\S3, pp.~32--35]{OtalPeigne2004}; see also the zero-potential case of \cite[Lemma~6.6(1)--(2)]{PPS2015}. We check it for the present kernels without assuming $\log\rho\in L^1(\nu)$.

For all integers $n\geq1$, define on a common full-measure set
\begin{displaymath}
    q_n(v)=\frac{m_{W^u(v)}^u(\zeta_n(v))}{\rho(v)}.
\end{displaymath}
The numerator is a Borel kernel integral as in \eqref{eq:rigidity-kernel-explicit}, using the code of $f^{-n}\zeta$, so $q_n$ is measurable. The actual set identity $f^n\zeta_n(v)=\zeta(f^nv)$ and \eqref{eq:rigidity-leafscale} give
\begin{displaymath}
    q_n(v)=e^{-\delta_{\Gamma} n\tau}\frac{\rho(f^nv)}{\rho(v)}.
\end{displaymath}
Subordination and reference positivity give $q_n>0$, while inclusion of the sub-atom gives $q_n\leq1$. Hence
\begin{displaymath}
    -\log q_n(v)
    =\delta_{\Gamma} n\tau+\log\rho(v)-\log\rho(f^nv)\geq0.
\end{displaymath}
The negative part of $G-G\circ T$, with $G=\log\rho$ and $T=f^n$, is bounded by $\delta_{\Gamma} n\tau$. Since $T$ preserves the probability $\nu$ and $G$ is finite almost everywhere, \cite[Lemma~8, p.~33]{OtalPeigne2004} applies with its function $H=-G$, without requiring $G\in L^1(\nu)$. It gives $G-G\circ T\in L^1(\nu)$ and $\int(G-G\circ T)\dd\nu=0$, and consequently
\begin{equation}\label{eq:rigidity-crossentropy}
    \int-\log q_n\dd\nu=\delta_{\Gamma} n\tau<\infty.
\end{equation}

\subsection{Entropy of the finite BMS measure}\label{sec:finite}

Suppose that $0<\mass m<\infty$ and write $\bar m=m/\mass m$. By \eqref{eq:prelim:finite_BMS}, this probability is flow ergodic and gives full measure to $\Reg$. Apply the preceding construction with $\nu=\bar m$, choosing $\tau>0$ so that $f=\phi^\tau$ is ergodic. We adapt the zero-potential disintegration argument of \cite[Lemma~6.5]{PPS2015} to the complete Borel atoms of $\zeta$.

\begin{lemma}\label{lem:finite-bms-conditionals}
    For the partition $\zeta$ constructed with $\nu=\bar m$, the conditional probabilities of $\bar m$ satisfy
    \begin{displaymath}
        \bar m_A=m_A^0=\frac{m^u|_A}{m^u(A)}
    \end{displaymath}
    for almost every $\zeta$ atom $A$.
\end{lemma}
\begin{proof}
In Hopf coordinates on the regular set, write $(\xi,\eta,t)=(v^-,v^+,\beta_{v^-}(\pi v,o))$. Since $G_o(\xi,\eta)=-t+\beta_\eta(o,\pi v)$, equations \eqref{eq:hopf_bms} and \eqref{eq:rigidity-leafmeasure} give
\begin{displaymath}
    \dd\widetilde m=e^{-\delta_{\Gamma} t}\dd\sigma_o(\xi)\dd t\,\dd m^u_{\xi,t}.
\end{displaymath}
Thus, local unstable plaque conditionals are normalized restrictions of the geometric leaf measure. The same holds in the transported foliation box $f^kW$, by flow invariance and \eqref{eq:rigidity-leafscale}.

Let $S_k=\{v:N(v)=k\}$. These Borel sets are saturated by $\zeta$, since $N(v)$ is the index of the first non-$*$ entry of the code $c(v)$, and their union has full $\bar m$ measure. By \eqref{eq:rigidity-plaque-container}, every atom in $S_k$ lies in a single plaque of $f^kW$. Its transverse parameter is the $k$th code entry, so the plaque partition is a measurable coarsening of $\zeta$ on $S_k$.

Fix $k$ with $\bar m(S_k)>0$ and disintegrate the normalized restriction of $\bar m$ to $S_k$ over these plaques. For almost every plaque $P$, its conditional probability is
\begin{displaymath}
    \lambda_{k,P}=\frac{m^u|_{P\cap S_k}}{m^u(P\cap S_k)}.
\end{displaymath}
The denominator is finite by bounded weights on $W$ and leaf scaling; parameters where it vanishes have zero transverse marginal.

By \eqref{eq:rigidity-rho} and disintegration, for almost every such plaque and $\lambda_{k,P}$-almost every $v$, the complete atom $A=\zeta(v)$ has positive finite leaf measure. Since $A\subset P\cap S_k$, it has positive $\lambda_{k,P}$ measure. There are at most countably many disjoint atoms of positive mass, and they cover a set of full $\lambda_{k,P}$ measure. Conditional probability on each is therefore
\begin{displaymath}
    \frac{\lambda_{k,P}|_A}{\lambda_{k,P}(A)}=\frac{m^u|_A}{m^u(A)}.
\end{displaymath}
The right-hand side is the measurable kernel from \eqref{eq:rigidity-kernel-explicit}. Transitivity and uniqueness of disintegration identify it with the conditional probability of $\bar m$ on $A$. Saturation of $S_k$ ensures that restricting to $S_k$ leaves these complete atoms unchanged. Combining the countably many sets $S_k$ proves the assertion.
\end{proof}

\begin{corollary}\label{cor:finite-bms-entropy}
    If $0<\mass m<\infty$, then
    \begin{displaymath}
        h_{\bar m}(\phi^1)=\delta_{\Gamma}.
    \end{displaymath}
\end{corollary}
\begin{proof}
By Lemma~\ref{lem:finite-bms-conditionals},
\begin{displaymath}
    q_1(v)=\bar m_{\zeta(v)}\bigl(\zeta_1(v)\bigr)
    \qquad\text{for }\bar m\text{-almost every }v.
\end{displaymath}
Since $q_1>0$ almost everywhere, the induced distribution on child atoms is countably supported with total mass one for almost every parent atom. Proposition~\ref{prop:rigidity-fullentropy} and \eqref{eq:rigidity-crossentropy} now give
\begin{displaymath}
    h_{\bar m}(f)=H_{\bar m}(f^{-1}\zeta\mid\zeta)=\int-\log q_1\dd\bar m=\delta_{\Gamma}\tau.
\end{displaymath}
The time-scaling formula \eqref{et:time_scaling} yields the result.
\end{proof}

\subsection{Entropy rigidity for regular probabilities}\label{sec:rigidity}

\begin{theorem}\label{thm:regular-rigidity}
    Under the standing hypotheses, if $\nu$ is a flow-ergodic invariant probability measure and $\nu(\Reg)=1$, then
    \begin{displaymath}
        h_\nu(\phi^1)\leq\delta_{\Gamma}.
    \end{displaymath}
    Equality implies
    \begin{displaymath}
        0<\mass m<\infty,
        \qquad \nu=\frac{m}{\mass m}.
    \end{displaymath}
\end{theorem}

We now fix an arbitrary flow-ergodic probability $\nu$ with $\nu(\Reg)=1$ and use the sampling time, partition, and reference kernels constructed for $\nu$ in Section~\ref{sec:unstable-tools}.

Disintegrate $\nu$ over $\zeta$ as $\nu=\int\nu_A\dd\bar\nu_\zeta(A)$. For a typical parent atom $A$ and $n\geq1$, subordination at integer translates gives leafwise neighborhoods in $\zeta_n(v)$ for $\nu_A$-almost every $v$. Disjointness of their nonempty interiors in a second-countable leaf makes these sub-atoms a countable family covering almost all of $A$. List the positive-mass complete Borel sub-atoms as $B_1,B_2,\ldots$, and put
\begin{displaymath}
    p_i=\nu_A(B_i),\qquad q_i=m_A^0(B_i),
    \qquad B_0=A\setminus\bigcup_{i\geq1}B_i.
\end{displaymath}
Then $p_0=0$, while the reference mass $q_0=m_A^0(B_0)$ may be positive. In particular,
\begin{displaymath}
    \sum_{i\geq1}p_i=1,
    \qquad \sum_{i\geq1}q_i\leq1.
\end{displaymath}
Each $B_i$ contains a typical center with positive endpoint in $\Lambda_\Gamma$, so subordination gives $q_i>0$.

Define
\begin{displaymath}
    p_n(v):=\nu_{\zeta(v)}(\zeta_n(v)).
\end{displaymath}
Integrating the Borel relation defining $\zeta_n$ against the conditional kernel makes $p_n$ measurable; on $B_i$, the functions $p_n$ and $q_n$ equal $p_i$ and $q_i$. By \eqref{eq:rigidity-crossentropy}, Proposition~\ref{prop:rigidity-fullentropy}, and the upper bound in Proposition~\ref{prop:up:entropy}, the conditional integrals of $-\log q_n$ and $-\log p_n$ are finite almost everywhere. Their difference is measurable, and Jensen's inequality \cite[\S3, p.~35]{OtalPeigne2004} gives
\begin{displaymath}
    D_n(A):=\sum_{i\geq1}p_i\log\frac{p_i}{q_i}
    \geq-\log\left(\sum_{i\geq1}q_i\right)\geq0.
\end{displaymath}
Equality holds precisely when $q_i=p_i$ for every $i\geq1$ and $q_0=0$. Integrating and using time scaling yields
\begin{equation}\label{eq:rigidity-KL}
    \int D_n(A)\dd\bar\nu_\zeta(A)=\delta_{\Gamma} n\tau-H_\nu(\zeta_n\mid\zeta)=n\tau\left(\delta_{\Gamma}-h_\nu(\phi^1)\right).
\end{equation}

Suppose now that $h_\nu(\phi^1)=\delta_{\Gamma}$. Equation~\eqref{eq:rigidity-KL} implies that, for almost every parent atom $A$, the equality $D_n(A)=0$ holds simultaneously for every $n\geq1$. Hence
\begin{equation}\label{eq:rigidity-atomequality}
    p_i=q_i\quad(i\geq1),\qquad q_0=0.
\end{equation}
For a bounded Lipschitz function $F$ on $Y$, put $\operatorname{osc}_B F=\sup_{x,y\in B}|F(x)-F(y)|$; since \eqref{eq:rigidity-atomequality} gives equal masses on the complete child atoms and zero mass to the remainder,
\begin{displaymath}
    \left|\int F\dd\nu_A-\int F\dd m_A^0\right|
    \leq\int\operatorname{osc}_{\zeta_n(v)}F\dd\nu_A(v)
    \longrightarrow0.
\end{displaymath}
The integrand is measurable on the countable full-measure union of the $B_i$, bounded by $2\|F\|_\infty$, and tends to zero by Lemma~\ref{lem:rigidity-subordinate}, so dominated convergence applies. A countable determining family of bounded Lipschitz functions therefore gives, as in \cite[\S3, p.~35]{OtalPeigne2004} and \cite[Lemma~6.6(5)]{PPS2015}, for almost every $A$,
\begin{equation}\label{eq:rigidity-conditionalequality}
    \nu_A=m_A^0.
\end{equation}

To pass to local leaf conditionals using \cite[\S6.1]{PPS2015}, we must allow a complete atom to meet a box in several local leaf pieces.

\begin{lemma}\label{lem:rigidity-localization}
    If $\nu_A\ll m_A^0$ for almost every $\zeta$ atom, then in every sufficiently small regular Hopf box the local strong unstable conditional measures of $\nu$ are absolutely continuous with respect to the geometric Patterson--Sullivan leaf measures.
\end{lemma}
\begin{proof}
Fix a relatively compact Hopf box $\mathcal B$, identified in Hopf
coordinates with $I_-\times I_+\times J$, with $\nu(\mathcal B)>0$. Let $q(\xi,\eta,t)=(\xi,t)$ be the transverse projection, and write $P_a=q^{-1}(a)$ for its unstable plaques.

By \eqref{eq:rigidity-plaque-container}, almost every atom $A$ lies in an embedded leaf disk $D$. The open subset $D\cap\mathcal B$ has at most countably many connected components, each contained in a single plaque. By consistency of the leaf measures, a Borel set $E\subset\mathcal B$ that is reference-null on every plaque satisfies $m_A^0(E)=0$; the hypothesis and disintegration then give $\nu(E)=0$.

Set $\nu_{\mathcal B}=\nu|_{\mathcal B}/\nu(\mathcal B)$ and $\alpha=q_*\nu_{\mathcal B}$. Typical positive endpoints lie in $\Lambda_\Gamma$, so $\sigma_o(I_+)>0$. The positive, locally bounded weights ensure that each reference plaque measure is finite and positive; the normalization gives a measurable probability kernel $\kappa_a$. Define
\begin{displaymath}
    \lambda(E)=\int\kappa_a(E)\dd\alpha(a).
\end{displaymath}

If $\lambda(E)=0$, then $N=\{a:\kappa_a(E)>0\}$ is Borel and $\alpha$-null. The set $E\setminus q^{-1}(N)$ is reference-null on every plaque, hence $\nu_{\mathcal B}$-null by the preceding paragraph. Since $\nu_{\mathcal B}(q^{-1}(N))=\alpha(N)=0$, this proves $\nu_{\mathcal B}\ll\lambda$.

Put $h=d\nu_{\mathcal B}/d\lambda$. The common transverse marginal $\alpha$ gives $\int h\dd\kappa_a=1$ almost everywhere, so $h\kappa_a$ disintegrates $\nu_{\mathcal B}$. Uniqueness gives $\nu_{\mathcal B,a}=h\kappa_a\ll\kappa_a$ almost everywhere.
\end{proof}

Apply the construction to $\phi^{-t}$: the critical exponent is unchanged, $h_\nu(f^{-1})=h_\nu(f)$, $f^{-1}$ is ergodic, and a positive-definite window can be chosen anew. Its strong unstable leaves are the original strong stable leaves, with reference measures
\begin{displaymath}
    \dd m^s(w)=e^{\delta_{\Gamma}\beta_{w^-}(o,\pi w)}\dd\sigma_o(w^-).
\end{displaymath}
They satisfy the same positive-exponent image-set scaling under the reversed flow. Entropy equality gives local strong stable conditional absolute continuity as well.

We deduce $\nu\ll m$ from the following local form of the flow-current correspondence \cite[Theorem~2.1]{Kaimanovich1990} and Fubini's theorem; only local finiteness is needed.

Lift $\nu$ to the locally finite, $\Gamma$-invariant and flow-invariant measure $\widetilde\nu$ on $T^1X$. On the regular Hopf domain, use
\begin{displaymath}
    t=\beta_\xi(\pi v,o),\qquad
    s=\beta_\eta(o,\pi v),\qquad
    (\xi,\eta)=(v^-,v^+).
\end{displaymath}
Both time coordinates increase at unit speed under the forward flow. They satisfy
\begin{equation}\label{eq:rigidity-timechange}
    s=t+G_o(\xi,\eta),\qquad
    G_o(\xi,\eta)=\beta_\xi(o,\pi v)+\beta_\eta(o,\pi v),
\end{equation}
where $G_o$ is independent of the point $v$ on the geodesic.

By the Hopf-coordinate identification~\eqref{eq:rigidity-hopfmap}, the regular Hopf image is $\mathcal E_{\Reg}\times\R$. For an endpoint set $D$ with compact closure in $\mathcal E_{\Reg}$, the inverse image of $D\times[0,1]$ has compact closure. Thus $I\mapsto\widetilde\nu(D\times I)$ is a locally finite translation-invariant measure on $\R$, equal to $\lambda(D)$ times Lebesgue measure. Taking $I=[0,1]$ defines a countably additive, locally finite endpoint measure $\lambda$. Rectangles and an exhaustion by relatively compact endpoint sets give
\begin{equation}\label{eq:rigidity-current}
    \dd\widetilde\nu=\dd\lambda(\xi,\eta)\dd t.
\end{equation}
The coordinate change \eqref{eq:rigidity-timechange} translates time separately on each endpoint fiber, so the same endpoint measure also gives $\dd\widetilde\nu=\dd\lambda\dd s$.

Restrict $\lambda$ to an endpoint rectangle $I_-\times I_+$ with closure in $\mathcal E_{\Reg}$, ignoring zero-mass rectangles, and normalize. Let $\alpha(d\xi)$ and $\beta(d\eta)$ be its marginals. Strong unstable plaques fix $(\xi,t)$, so their conditional absolute continuity and \eqref{eq:rigidity-current} imply
\begin{displaymath}
    \lambda(d\eta\mid\xi)\ll\sigma_o(d\eta)
    \quad\text{for }\alpha\text{-almost every }\xi.
\end{displaymath}
The first application of Fubini gives
\begin{equation}\label{eq:rigidity-fubini1}
    \lambda\ll\alpha\otimes\sigma_o.
\end{equation}
Here and below, boundary measures are restricted to the corresponding endpoint neighborhoods.

In the $s$ coordinates, strong stable plaques fix $(\eta,s)$. Their conditional absolute continuity similarly gives
\begin{displaymath}
    \lambda(d\xi\mid\eta)\ll\sigma_o(d\xi)
    \quad\text{for }\beta\text{-almost every }\eta.
\end{displaymath}
Integrating these conditionals against $\beta$ gives $\alpha(A_-)=0$ whenever $\sigma_o(A_-)=0$. Thus $\alpha\ll\sigma_o$, and \eqref{eq:rigidity-fubini1} yields
\begin{displaymath}
    \lambda\ll\sigma_o\otimes\sigma_o,
    \qquad
    \widetilde\nu\ll\sigma_o(d\xi)\sigma_o(d\eta)\dd t.
\end{displaymath}

By \eqref{eq:hopf_bms}, the BMS density relative to this product is positive and finite, so $\nu\ll m$ on every small regular box. A countable cover of $\Reg$ and $\nu(\Reg)=1$ gives global absolute continuity.

\begin{proof}[Completion of the proof of Theorem~\ref{thm:regular-rigidity}]
    Equation~\eqref{eq:rigidity-KL} gives $h_\nu(\phi^1)\leq\delta_{\Gamma}$. Under equality, \eqref{eq:rigidity-conditionalequality} and Lemma~\ref{lem:rigidity-localization}, applied in both time directions, give local conditional absolute continuity. The preceding product-measure argument yields $\nu\ll m$, and \eqref{eq:prelim:ac_BMS} gives
    \begin{displaymath}
        0<\mass m<\infty,\qquad \nu=m/\mass m.
    \end{displaymath}
\end{proof}

\subsection{Singular measures and sparse returns}\label{sec:sampling}
\label{sec:sparse_sampling}

We now prove that every invariant probability supported on $\Sing$
has entropy strictly less than $\delta_{\Gamma}$ when $m$ is finite. Note that singular here means concentrated on $\Sing$, not singular with respect to $m$. Finiteness implies $\Gamma$ is divergent by \eqref{eq:prelim:finite_BMS}. The return coding and entropy reduction require neither assumption and will also be used in Theorem~\ref{thm:main}.

The entropy-deficit comparison parallels the compact uniqueness argument presented in \cite[Section~8]{Ricks2021}, following Knieper \cite{Knieper1998}. The return coding on a fixed compact set replaces global partitions, while conservative recurrence and singular-orbit avoidance give a vanishing factor uniform in the return step.
\begin{theorem}\label{thm:np:exclusion}
    Under the standing hypotheses, suppose that $\mass m<\infty$. Then every probability measure $\nu$ invariant under the geodesic flow with $\nu(\Sing)=1$ satisfies
    \begin{displaymath}
        h_\nu(\phi^1)<\delta_{\Gamma}.
    \end{displaymath}
\end{theorem}

Fix a rank one axis and choose $o$ on it. We first treat flow-ergodic probabilities. Keeping a positive-measure compact return set fixed while increasing the sampling time gives a return-label bound of the form $C\varepsilon_T e^{\delta_{\Gamma}nT}$, with $\varepsilon_T\to0$ as $T\to\infty$

Kac's and Abramov's formulas give the return-map normalization and entropy scaling (Section~\ref{sec:entropy_tools}; \cite{AbramovInd59,Abramov59}). The endpoint covering estimate shows that the partition defined by the return time and the corresponding deck transformation carries the full entropy of the induced map, even when the return set is a Borel subset of $\Sing$, without assuming that this partition is generating.

\begin{lemma}\label{lem:np:return_entropy}
    Suppose that $\tau>0$ and that $f=\phi^\tau$ preserves a probability measure $\nu$ and is ergodic. Let $A$ be a Borel subset of a relatively compact, evenly covered chart for $T^1X\to T^1M$, with $p=\nu(A)>0$. Fix the lift $L:A\to T^1X$ determined by that chart, and assume that $\overline{L(A)}$ is compact. Set
    \begin{displaymath}
        N(v)=\min\{n\ge1:f^nv\in A\},\qquad
        F(v)=f^{N(v)}v,\qquad \mu=\nu|_A/p.
    \end{displaymath}
    On the full-measure domain of finite returns, let $\gamma(v)\in\Gamma$ be the unique element satisfying
    \begin{displaymath}
        \phi^{\tau N(v)}L(v)=\gamma(v)L(Fv),
    \end{displaymath}
    and let $\cP$ be the countable partition with label $(N,\gamma)$. If $\delta_{\Gamma}<\infty$, then
    \begin{equation}\label{eq:np:return_entropy}
        H_\mu(\cP)<\infty,
        \qquad h_\mu(F)=h_\mu(F,\cP)\le H_\mu(\cP).
    \end{equation}
    If $\nu$ is invariant under the full flow, then
    \begin{equation}\label{eq:np:abramov}
        \mathbb E_\mu N=\frac1p,
        \qquad \tau h_\nu(\phi^1)=p\,h_\mu(F).
    \end{equation}
\end{lemma}

\begin{proof}
    By recurrence, $N<\infty$ almost everywhere; on the recurrent full-measure domain, $F$ is an invertible, measure-preserving, ergodic transformation of $(A,\mu)$.

    Let $E$ be the standard projection of $\overline{L(A)}$, and choose $D$ with $E\subset B(o,D)$. If $N=n$, the lifted orbit segment has length $\tau n$ and joins $E$ to $\gamma E$. Consequently,
    \begin{displaymath}
        d_X(o,\gamma o)\le\tau n+2D.
    \end{displaymath}
    For each $s>\delta_{\Gamma}$, the definition of the critical exponent supplies $C_s\ge1$ such that the number of possible deck labels at step $n$ is at most $C_s e^{s\tau n}$.

    Kac's formula~\eqref{et:kac} gives $\mathbb E_\mu N=p^{-1}$. Thus, the label estimate~\eqref{et:finite_labels} yields
    \begin{displaymath}
        H_\mu(\cP)
        \le\log C_s+(\log2+s\tau)\mathbb E_\mu N<\infty.
    \end{displaymath}

    Write $\cP^k=\bigvee_{j=0}^{k-1}F^{-j}\cP$. On each atom, the cumulative return steps and deck transformations
    \begin{displaymath}
        S_j=\sum_{i=0}^{j-1}N\circ F^i,
        \qquad
        \Gamma_j=\gamma\, (\gamma\circ F)\cdots(\gamma\circ F^{j-1})
        \quad(0\le j\le k)
    \end{displaymath}
    are fixed, with $S_0=0$ and $\Gamma_0=e$. Since the flow commutes with deck transformations, successive substitution gives
    \begin{displaymath}
        \phi^{\tau S_j}L(v)=\Gamma_j L(F^jv).
    \end{displaymath}
    All lifted segments therefore join $E$ to $\Gamma_kE$ in time $\tau S_k$. Lemma~\ref{lem:up:endpoint} divides the atom into at most $C_\eps=C(E,E,\eps/2)$ sets of lifted Bowen diameter less than $\eps$ on $[0,\tau S_k]$, independently of $k,S_k,\Gamma_k$. Projection at the common observation times $\tau S_j$ makes their diameter less than $\eps$ also for
    \begin{displaymath}
        d_k^F(v,w)=\max_{0\le j<k}d(F^jv,F^jw),
    \end{displaymath}
    where $d$ is the quotient of the geodesic metric in \eqref{eq:up:geodesicmetric}.

    We use the standard compact-core conditional-entropy comparison; see \cite[proof of Theorem~3.4.1]{PU2010}. The names are estimated by Hamming counting within each return-code atom. Fix a finite partition $\cQ=\{Q_1,\ldots,Q_q\}$, $q\geq2$, write $\cQ^k=\bigvee_{j=0}^{k-1}F^{-j}\cQ$, and fix $0<\kappa<1/16$. Choose compact cores $K_i\subset Q_i$ with $\mu(K)>1-\kappa$, where $K=\bigcup_iK_i$, and choose $\eps>0$ so that distinct nonempty cores have distance greater than $4\eps$. By ergodicity of $F$, the sets
\begin{displaymath}
 G_k=\left\{v:\sum_{j=0}^{k-1}\indic_{A\setminus K}(F^jv)
                     \leq2\kappa k\right\}
\end{displaymath}
satisfy $\mu(G_k)\to1$. In each of the at most $C_\eps$ Bowen-small sets covering a fixed $\cP^k$ atom, two names from $G_k$ differ at at most $4\kappa k$ positions. Thus, the number of such $\cQ^k$-names in that atom is at most
\begin{displaymath}
 C_\eps\sum_{j=0}^{\lfloor4\kappa k\rfloor}\binom{k}{j}q^j
 \leq C_\eps(k+1)
       \exp\{k[b(4\kappa)+4\kappa\log q]\},
\end{displaymath}
where $b(u)=-u\log u-(1-u)\log(1-u)$. Recording membership in $G_k$ costs at most $\log2$; there are at most $q^k$ names on its complement. Therefore
\begin{displaymath}
 H_\mu(\cQ^k\mid\cP^k)
 \leq\log2+\log C_\eps+\log(k+1)
      +k[b(4\kappa)+4\kappa\log q]
      +\mu(G_k^c)k\log q.
\end{displaymath}
Using $H_\mu(\cQ^k)\leq H_\mu(\cP^k)+H_\mu(\cQ^k\mid\cP^k)$, first let $k\to\infty$ and then $\kappa\downarrow0$ to obtain $h_\mu(F,\cQ)\leq h_\mu(F,\cP)$. Taking the supremum over finite $\cQ$, and using the reverse standard inequality for countable finite-entropy partitions, proves \eqref{eq:np:return_entropy}. Kac's formula, \eqref{et:abramov_induced}, and \eqref{et:time_scaling} give \eqref{eq:np:abramov}.
\end{proof}

For a full-flow ergodic invariant probability $\nu$, a countable cover by relatively compact evenly covered charts and inner regularity give a compact $A$ in one chart with $p=\nu(A)>0$. If $\nu(\Sing)=1$, choose $A\subset\Sing$. Fix its lift $L(A)$ and compact standard projection $E\subset B(o,D)$. Keep $A,p,E,D$ fixed as the sampling time varies.

By the Pugh--Shub theorem \cite[Theorem~1]{PughShub71}, as recalled in Section~\ref{sec:entropy_tools}, there are arbitrarily large $T>0$ for which $\phi^T$ is ergodic. For each such $T$, let $N_T,F_T,\cP_T$ be the return-step function, induced map, and label partition on the same set $A$, and put
\begin{displaymath}
    \mu=\nu|_A/p,
    \qquad r_{n,T}=\mu(N_T=n).
\end{displaymath}
Equations~\eqref{et:kac} and \eqref{et:finite_labels} give the uniform bound
\begin{displaymath}
    pH(r_{\cdot,T})
    \le p(\log2)\mathbb E_\mu N_T=\log2.
\end{displaymath}

Let $\Gamma_{n,T}$ consist of the deck transformations whose label $(n,\gamma)$ has positive $\mu$-mass, and put $M_{n,T}=\#\Gamma_{n,T}$, finite and positive whenever $r_{n,T}>0$. A realizing point supplies a minimizing lifted segment of length $nT$ from $E$ to $\gamma E$, so
\begin{equation}\label{eq:np:annulus}
    \left|d_X(o,\gamma o)-nT\right|\le2D.
\end{equation}
The entropy chain rule and Lemma~\ref{lem:np:return_entropy} imply
\begin{align}
    T h_\nu(\phi^1)
 &=p h_\mu(F_T)\le pH_\mu(\cP_T)\notag\\
 &\le pH(r_{\cdot,T})
 +p\sum_{n:r_{n,T}>0}r_{n,T}\log M_{n,T}.
 \label{eq:np:one_label_entropy}
\end{align}
It remains to find $C>0$, independent of $n,T$, and $\eps_T\to0$ such that, for all sufficiently large ergodic sampling times,
\begin{equation}\label{eq:np:count_target}
    M_{n,T}\le C\eps_T e^{\delta_{\Gamma} nT}
    \qquad(n\ge1).
\end{equation}
Recurrence to $A$ rules out $\eps_T=0$, which would eliminate every positive-mass return label. Substituting \eqref{eq:np:count_target} into \eqref{eq:np:one_label_entropy} gives
\begin{align}
    T h_\nu(\phi^1)
 &\le\log2+p\log(C\eps_T)
 +p\delta_{\Gamma} T\sum_{n\ge1}nr_{n,T}\notag\\
 &=\delta_{\Gamma} T+\log2+p\log(C\eps_T).
 \label{eq:np:master}
\end{align}
Here $p\mathbb E_\mu N_T=1$. Since $p>0$ is fixed, $\eps_T\to0$ makes the sum of the last two terms strictly negative for large sampling times.

For the geometric estimate, write $\pi_X:T^1X\to X$ and $\pi_\Gamma:T^1X\to T^1M$ for the standard and covering projections.

The required tube bound follows directly from the individual-element overlap estimate in \cite[Section~6, proof of Proposition~6, before the summation giving~(32)]{LP16}. With $o$ on the fixed rank one axis, that argument provides $c\ge1$, $a>0$, $R_0>4c$, and a compact set $\mathcal K_0\subset\pi_X^{-1}(\overline B(o,2c))$ such that
\begin{displaymath}
    \int_0^\infty\widetilde m(\mathcal K_0\cap\phi^{-t}\gamma\mathcal K_0)\,\dd t\ge a e^{-\delta_{\Gamma} d_\gamma},\qquad d_\gamma=d_X(o,\gamma o)>R_0.
\end{displaymath}
Here, $G_o=2\operatorname{Gr}_o$, so its geometric measure is exactly \eqref{eq:full_knieper}. The estimate uses the rank one endpoint neighborhoods, finite-correction corridor, and shadow bound and requires neither divergence of the Poincar\'{e} series at $\delta_{\Gamma}$ nor finiteness of the quotient measure.

Set $\mathcal K=\pi_X^{-1}(\overline B(o,2c))$, which is compact by Hopf--Rinow and compactness of the unit sphere fibers. Enlarging $\mathcal K_0$ to $\mathcal K$ gives
\begin{equation}\label{eq:np:integrated_tube}
    \int_0^\infty\widetilde m(\mathcal K\cap\phi^{-t}\gamma\mathcal K)\,\dd t
    \ge a e^{-\delta_{\Gamma} d_\gamma}.
\end{equation}
A nonempty overlap at time $t\ge0$ supplies a minimizing geodesic segment joining the two closed balls, and hence
\begin{displaymath}
    |t-d_\gamma|\le4c.
\end{displaymath}
Define
\begin{displaymath}
    D_\gamma=\{w\in\mathcal K:\phi^tw\in\gamma\mathcal K~\text{for some}~t\in[d_\gamma-4c,d_\gamma+4c]\}.
\end{displaymath}
It is compact, being the projection of the closed set $\{(w,t):\phi^tw\in\gamma\mathcal K\}$ in $\mathcal K\times[d_\gamma-4c,d_\gamma+4c]$. The overlaps in \eqref{eq:np:integrated_tube} lie in $D_\gamma$ and have time support of length $8c$, so
\begin{equation}\label{eq:np:tube_lower}
    \widetilde m(D_\gamma)\ge\frac{a}{8c}e^{-\delta_{\Gamma} d_\gamma}
    \qquad(d_\gamma>R_0).
\end{equation}
The constants depend only on the fixed geometric data and the Patterson-Sullivan density, not on the candidate invariant probability.

The finite-window estimate below parallels the separation argument in \cite[Section~8]{Ricks2021}, using compactness and the flat strip theorem \cite[Corollary~I.5.8(ii)]{Ballmann95}. Recall that the singular sets are closed and invariant under the flow and time reversal, and the lifted set is invariant under deck transformations.

Choose a vector $v_0$ on the given closed rank one geodesic and an open neighborhood $U$ with compact closure in a neighborhood in $\Reg$ equipped with Hopf coordinates, using the identification \eqref{eq:rigidity-hopfmap}. Its lift contains an open endpoint rectangle about the axis endpoints times a nondegenerate time interval. Since those endpoints belong to $\Lambda_\Gamma=\supp\sigma_o$, the Hopf formula~\eqref{eq:hopf_bms} gives $m(U)>0$. Local finiteness then gives
\begin{displaymath}
    v_0\in U,\qquad \overline U\subset\Reg,
    \qquad 0<m(U)<\infty.
\end{displaymath}

\begin{lemma}\label{lem:np:avoidance}
    For every $R<\infty$, there is $L=L(U,R)>0$ with the following property. If $c$ is a complete singular geodesic and $w$ is another complete unit-speed geodesic satisfying
    \begin{displaymath}
        d_X(c(t),w(t))\le R\qquad(0\le t\le q),
    \end{displaymath}
    then $\pi_\Gamma(w'(t))\notin U$ for every $L\le t\le q-L$.
\end{lemma}

\begin{proof}
If no uniform $L$ existed, there would be singular geodesics $c_j$, tracking geodesics $w_j$, and times $t_j\in[0,q_j]$ with $\min\{t_j,q_j-t_j\}\to\infty$ and $\pi_\Gamma(w_j'(t_j))\in U$. Shift time by $t_j$ and apply common deck transformations so that $w_j'(0)$ lies in a fixed finite union of compact lifts covering $\overline U$. The tracking bound places the base points of $c_j'(0)$ in a fixed compact neighborhood. The properness and compactness of unit sphere fibers implies convergent subsequences of both initial vectors.

The complete limiting geodesics $c_\infty,w_\infty$ satisfy $d_X(c_\infty(t),w_\infty(t))\leq R$ for every $t\in\R$. Closedness of $\Sing$ makes $c_\infty$ singular, whereas $\pi_\Gamma(w_\infty'(0))\in\overline U\subset\Reg$. The flat strip theorem \cite[Corollary~I.5.8(ii)]{Ballmann95} and regularity of $w_\infty$ make the two limiting geodesics identical up to time shift, contradicting flow invariance of rank.
\end{proof}

Assume now that the Poincar\'{e} series diverges at $\delta_{\Gamma}$ and $\nu(\Sing)=1$, with $A\subset\Sing$ and $p,E,D$ fixed above. For each $\gamma\in\Gamma_{n,T}$, let $c_\gamma$ be the complete singular lift of a point realizing the label. With $q=nT$,
\begin{displaymath}
    c_\gamma(0)\in E,\qquad c_\gamma(q)\in\gamma E.
\end{displaymath}
Set $R=6c+3D$ and $L=L(U,R)$, and restrict to ergodic sampling times with
\begin{equation}\label{eq:np:T_threshold}
    T>\max\{2L,\ R_0+2D,\ 4c+2D\}.
\end{equation}
All counted displacements then exceed $R_0$.

Take $w\in D_\gamma$. There is $s>0$ such that $\phi^sw\in\gamma\mathcal K$ and $|s-d_\gamma|\le4c$. Together with \eqref{eq:np:annulus}, this gives $|s-q|\le4c+2D$. For the geodesic $c_w$ determined by $w$, we obtain
\begin{align*}
    d_X(c_w(0),c_\gamma(0))&\le2c+D,\\
    d_X(c_w(q),c_\gamma(q))
                           &\le |q-s|+2c+D\le6c+3D.
\end{align*}
Convexity of distance between equally parametrized geodesics \cite[Proposition~I.5.4]{Ballmann95} implies
\begin{equation}\label{eq:np:coarse_tracking}
    d_X(c_w(t),c_\gamma(t))\le R\qquad(0\le t\le q).
\end{equation}
The same calculation gives
\begin{displaymath}
    d_X(c_w(q),\gamma o)\le |q-s|+2c\le6c+2D.
\end{displaymath}
The orbit-counting bound used in the proof of Lemma~\ref{lem:semigroup_entropy} shows that, for fixed $n,T,w$, at most $J$ elements $\gamma\in\Gamma_{n,T}$ satisfy $w\in D_\gamma$, where
\begin{displaymath}
    J=N_o(12c+4D)<\infty.
\end{displaymath}
Properness implies finiteness, independently of $n,T$. Combining this with \eqref{eq:np:tube_lower} and the annulus estimate yields
\begin{align}
    M_{n,T}\frac{a}{8c}e^{-2\delta_{\Gamma} D}e^{-\delta_{\Gamma} nT}
 &\le\sum_{\gamma\in\Gamma_{n,T}}\widetilde m(D_\gamma)\notag\\
 &\le J\,\widetilde m\left(
 \bigcup_{\gamma\in\Gamma_{n,T}}D_\gamma\right).
 \label{eq:np:tube_count}
\end{align}
All these sets lie in the same lifted compact set $\mathcal K$.

The conservative ergodic consequence \eqref{eq:prelim:HTS_visits} implies that almost every orbit visits $U$ at arbitrarily late positive times. This full-measure statement lifts to $\mathcal K$ through a countable atlas of evenly covered charts.

For $T>2L$, define
\begin{displaymath}
    \mathcal E_T=
    \{w\in\mathcal K:\pi_\Gamma(\phi^tw)\notin U
    \text{ for every }t\in[L,T-L]\}.
\end{displaymath}
Openness of $U$ and continuity of the flow make $\mathcal E_T$ closed in $\mathcal K$, hence compact and Borel. These sets decrease with $T$. Their intersection consists of vectors with no visit to $U$ after time $L$, and therefore has zero $\widetilde m$-measure.

\begin{figure}[htbp]
    \centering
    \begin{tikzpicture}[x=.95cm,y=.8cm,>=Stealth,font=\small]
        \draw[->,black!55] (0,0)--(8.6,0) node[right] {$t$};
        \foreach \x/\lab in {0/0,1/L,3/{T-L},4/T,7/{nT-L},8/nT}
        {\draw[black!55] (\x,.08)--(\x,-.08);
        \node[below=3pt] at (\x,0) {$\lab$};}
        \draw[surfaceedge,line width=2pt] (0,2.15)--(8,2.15);
        \node[above=4pt] at (4,2.15)
        {$w\in D_\gamma$: $R$-tracking of $c_\gamma$};
        \draw[glueboundary,line width=2pt] (1,1.15)--(7,1.15);
        \fill[glueboundary] (1,1.15) circle (1.7pt);
        \fill[glueboundary] (7,1.15) circle (1.7pt);
        \node[above=4pt] at (4,1.15) {avoidance of $U$};
        \draw[surfaceedge,line width=4pt] (1,.35)--(3,.35);
        \draw[densely dashed,black!35] (1,0)--(1,2.15);
        \draw[densely dashed,black!35] (3,0)--(3,.65);
        \draw[densely dashed,black!35] (7,0)--(7,2.15);
        \node[anchor=west] at (3.3,.4) {fixed interval for $\mathcal E_T$};
    \end{tikzpicture}
    \caption{The common avoidance interval for return times $nT$.}\label{fig:singular:avoidance}
\end{figure}
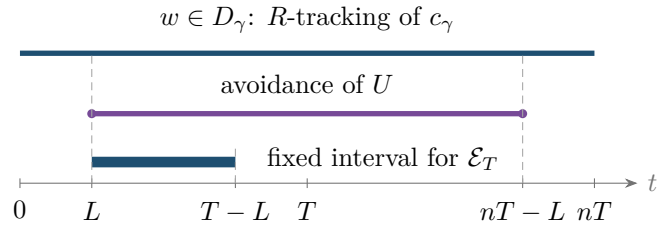

Local finiteness, as established in \eqref{eq:prelim:local_finiteness}, gives
\begin{displaymath}
    \widetilde m(\mathcal K)<\infty.
\end{displaymath}
Consequently,
\begin{displaymath}
    \eps_T:=\widetilde m(\mathcal E_T)\longrightarrow0.
\end{displaymath}
Since $q=nT\ge T$, \eqref{eq:np:coarse_tracking} and Lemma~\ref{lem:np:avoidance} give avoidance of $U$ on $[L,q-L]\supset[L,T-L]$. Thus
\begin{displaymath}
    D_\gamma\subset\mathcal E_T
    \qquad(\gamma\in\Gamma_{n,T}).
\end{displaymath}
The common interval $[L,T-L]$ makes $\eps_T$ independent of the return step $n$ (Figure~\ref{fig:singular:avoidance}).

Substitution into \eqref{eq:np:tube_count} gives
\begin{equation}\label{eq:np:final_count}
    M_{n,T}\le C\eps_T e^{\delta_{\Gamma}nT},
    \qquad C=\frac{8cJ}{a}e^{2\delta_{\Gamma} D}.
\end{equation}
The constant depends on the fixed return set, group, and density but not on $n,T$, proving \eqref{eq:np:count_target}.

Choose an ergodic sampling time satisfying \eqref{eq:np:T_threshold} and making $\log2+p\log(C\eps_T)<0$. Equations~\eqref{eq:np:final_count} and \eqref{eq:np:master} then give $h_\nu(\phi^1)<\delta_{\Gamma}$ for every singular full-flow ergodic probability.

The entropy-deficit identity \eqref{eq:et:entropy_deficit}, applied to the invariant Borel set $\Sing$, extends this strict bound to all invariant probabilities in the statement. This proves Theorem~\ref{thm:np:exclusion}. The entropy deficit for a singular invariant probability may depend on that probability.

Although stated for finite $m$, the argument uses only local finiteness and the conservative ergodic consequence \eqref{eq:prelim:HTS_visits}, and therefore gives the same strict bound whenever the Poincar\'{e} series diverges at $\delta_{\Gamma}$, even if $\mass m=\infty$.

\subsection{Proof of Theorem~\ref{thm:main}}
\begin{proof}
    Section~\ref{sec:entropy} gives
    \begin{displaymath}
        h_{\mathrm{top}}^{\,c}(\phi)=H=\delta_{\Gamma}.
    \end{displaymath}
    For any measure of maximal entropy $\nu$, write its full-flow ergodic decomposition as $\nu=\int\nu_\omega\dd\lambda(\omega)$. The entropy upper bound and \eqref{eq:et:entropy_equality} give $h_{\nu_\omega}(\phi^1)=\delta_{\Gamma}$ for almost every $\omega$. Since $\Reg$ and $\Sing$ are invariant, each ergodic component is concentrated on one of them.

    Suppose first that $\mass m<\infty$. Corollary~\ref{cor:finite-bms-entropy} gives $h_{\bar m}(\phi^1)=\delta_{\Gamma}$, so a measure of maximal entropy exists. Theorem~\ref{thm:np:exclusion} excludes singular maximizing components, and Theorem~\ref{thm:regular-rigidity} identifies every regular maximizing component with $\bar m$. Thus $\nu_\omega=\bar m$ for almost every $\omega$, and integration yields $\nu=\bar m$, proving uniqueness.

    Suppose now that $\mass m=\infty$. By Theorem~\ref{thm:regular-rigidity}, no regular ergodic probability can have entropy $\delta_{\Gamma}$. Hence, a measure of maximal entropy, if one existed, would have a singular ergodic component $\nu_\omega$ of entropy $\delta_{\Gamma}$. If the Poincar\'{e} series diverges at $\delta_{\Gamma}$, the observation at the end of Section~\ref{sec:sampling} excludes such a component even though $m$ is infinite.

    It remains to treat the case in which the Poincar\'e series converges at $\delta_{\Gamma}$. Use the return-coding setup of Section~\ref{sec:sampling}, with $A,p,E,D$ fixed for this component and $T$ ranging over its large ergodic sampling times. Set
    \begin{displaymath}
        \eps_T=e^{2\delta_{\Gamma} D}
        \sum_{d_X(o,\gamma o)\ge T-2D}e^{-\delta_{\Gamma} d_X(o,\gamma o)}.
    \end{displaymath}
    Convergence of the series gives $\eps_T\to0$. For $n\ge1$, every counted element satisfies $d_X(o,\gamma o)\ge nT-2D\ge T-2D$. Using both bounds in \eqref{eq:np:annulus}, we obtain
    \begin{displaymath}
        M_{n,T}e^{-\delta_{\Gamma} nT}
        \le e^{2\delta_{\Gamma} D}\sum_{\gamma\in\Gamma_{n,T}}
        e^{-\delta_{\Gamma} d_X(o,\gamma o)}
        \le\eps_T.
    \end{displaymath}
    This proves \eqref{eq:np:count_target} with $C=1$, uniformly in $n$. Equation~\eqref{eq:np:master}, at a sufficiently large ergodic sampling time, gives $h_{\nu_\omega}(\phi^1)<\delta_{\Gamma}$, a contradiction.

    Thus, no measure of maximal entropy exists when $\mass m=\infty$. Every invariant probability then has entropy strictly below $\delta_{\Gamma}$, while Theorem~\ref{thm:H_delta} supplies compactly supported invariant probabilities with entropies tending to $\delta_{\Gamma}$.
\end{proof}

\section{Mixing and the surface case}\label{sec:mixing}\label{sec:surfaces}

\begin{proof}[Proof of Theorem~\ref{thm:scope:mixing}]
    We apply the mixing theorem of Link~\cite[Theorem~A]{Link2020}, checking its geometric hypotheses and the identification of the invariant measure. The universal cover $X$ is a Hadamard manifold, which is proper and geodesically complete. The deck transformation group $\Gamma$ is discrete and non-elementary, and its limit set is infinite. Moreover, the assumed rank one axis cannot bound a flat half-plane, so the corresponding axial element is rank one in Link's terminology. Its fixed points belong to $\Lambda_\Gamma$. Hence, the rank one group hypothesis of Link's theorem is satisfied. The required finiteness of the critical exponent is part of our standing assumptions.

    The remaining dynamical hypothesis is divergence. By \eqref{eq:prelim:finite_BMS} and Theorem~\ref{thm:geo:HTS}, $\Gamma$ is divergent, and the flow is conservative and ergodic with respect to $m$.

    For the fixed Patterson--Sullivan density $(\sigma_x)_{x\in X}$, the construction in \eqref{eq:full_knieper} is precisely Knieper's construction. In the conservative manifold case, \cite[Theorem~C and the end of \S10]{Link2018} identifies it with Ricks's Bowen--Margulis measure constructed from the same density. The comparison takes place on a set of full measure consisting of zero-width geodesics where the parallel set is a single line. Its Riemannian volume is Lebesgue measure in the flow coordinates, as in \eqref{eq:hopf_bms}. Thus, the identification preserves the measure itself. The space of parametrized geodesic lines in $X$, modulo $\Gamma$, is naturally $T^1M$, and reparametrization is our flow $\phi^t$.

    Finally, the length spectrum in \cite{Link2020} consists of the translation lengths of axial elements. Translation along an invariant geodesic has length $\inf_{x\in X}d(x,\gamma x)$, so this is exactly $\mathcal L(\Gamma)$. The assumed non-arithmeticity rules out the second alternative in Link's theorem. We therefore obtain, for all Borel sets $A,B\subset T^1M$,
    \begin{displaymath}
        m(A\cap\phi^{-t}B)\longrightarrow\frac{m(A)m(B)}{\mass m}\qquad(t\to+\infty).
    \end{displaymath}
    All such sets have finite measure. Dividing by $\mass m$ gives the asserted limit for $\bar m=m/\mass m$.
\end{proof}

We now turn to the surface case. Let $S$ be a connected, orientable, smooth, complete surface without boundary, with $K\le0$ and genus at least two. The existence of a closed rank one geodesic and non-arithmeticity of the axial length spectrum follow from the geometric results of Link, Peign\'e, and Picaud~\cite{LPP2006,LP16}. We reduce to a finite-type covering surface, so no restriction on the topological type or volume of $S$ is needed.

Put $X=\widetilde S$ and $\Gamma=\pi_1(S)$. Choose a compact essential one-holed torus in $S$, and let $H\simeq F_2$ be its fundamental group, viewed as a subgroup of $\Gamma$. The covering surface $S_H=X/H$ is complete, orientable, noncompact, and of finite type in the sense of \cite{LPP2006}. It is not flat, since a free orientation-preserving action on the Euclidean plane is by translations and cannot contain $F_2$. By \cite[Lemma~1.18, Corollary~1.20, and Remark~1.21]{LPP2006}, $H$ has an infinite limit set, acts minimally on it, and contains an axial element $h$ whose axis bounds no flat half-plane. Any common boundary fixed point would be one of $h^-,h^+$ by \cite[Lemma~III.3.3]{Ballmann95}, contradicting minimality. Thus $H$, and hence $\Gamma$, is non-elementary; moreover, $\delta_\Gamma>0$ by \eqref{eq:prelim:free_growth}.

In dimension two, the existence of a weak rank one geodesic on a non-elementary quotient implies the existence of a closed rank one geodesic; see the discussion of the rank one hypotheses in \cite[Introduction]{LP16}. Applying this surface result to $S_H$ gives an axial element of $H$, and therefore of $\Gamma$, with a rank one axis. 

Finally, \cite[Lemma~3.2 and its proof, p.~35]{LPP2006} gives non-arithmeticity of $\mathcal L(H)$. Although that paper's length spectrum also includes exceptional parabolic displacements, the proof uses only axial elements. Since $H$ and $\Gamma$ act on the same $X$, their translation lengths agree on $H$, and $\mathcal L(H)\subset\mathcal L(\Gamma)$. Hence, $\mathcal L(\Gamma)$ is non-arithmetic as well.

\begin{proof}[Proof of Corollary~\ref{cor:scope:surface}]
The preceding geometric facts verify the hypotheses of Theorems~\ref{thm:main} and~\ref{thm:scope:mixing}; finiteness of $\delta_\Gamma$ is assumed. The two theorems give all the assertions of the corollary.
\end{proof}

\section*{Acknowledgments}
Fei Liu is partially supported by the National Natural Science Foundation of China under Grant No.~12571179.

Xiaokai Liu would like to thank Professor Raul Ures and Professor Jana Hertz for their invaluable guidance and continuous support at Southern University of Science and Technology.

\printbibliography
\end{document}